\documentclass[poms,final,nonblindrev]{poms1V1}

\OneAndAHalfSpacedXI

\usepackage{graphicx}
\usepackage{subfigure, epsfig}
\usepackage{natbib}

 \bibpunct[, ]{(}{)}{,}{a}{}{,}
 \def\bibfont{\small}

\TheoremsNumberedThrough
\ECRepeatTheorems

\EquationsNumberedThrough

\usepackage{hyperref}
    \hypersetup{
        colorlinks   = true,
        citecolor    = blue,
        linkcolor=blue
    }
\usepackage{adjustbox}
\usepackage{placeins}
\usepackage{array, makecell}
\usepackage{color, colortbl}
\definecolor{Gray}{gray}{0.9}
\usepackage{setspace,lipsum}
\usepackage{style}
\usepackage{xcolor}
\usepackage{arydshln}
\usepackage{pdflscape}
\usepackage{rotating}

\begin{document}

\RUNAUTHOR{Kurbanzade, Baig and Mehrotra}
\RUNTITLE{UAVs, Industries, Sustainability and Operational Innovation}
\TITLE{Elevating Industries with Unmanned Aerial Vehicles: Integrating Sustainability and Operational Innovation}

\ARTICLEAUTHORS{
\AUTHOR{Ali Kaan Kurbanzade {$^{a,c,d}$} $\bullet$ Ansaar M. Baig {$^{b,c,d}$} $\bullet$ Sanjay Mehrotra {$^{c,d,\ast}$}}
\AFF{$^a$ Department of Operations \& Decision Technologies, Kelley School of Business, Indiana University, Bloomington, IN, USA \\
$^b$ Harvard Business School, Harvard University, Boston, MA, USA \\
$^c$ Department of Industrial Engineering and Management Sciences, Northwestern University, Evanston, IL, USA \\
$^d$ Center for Engineering and Health, Northwestern University, Feinberg School of Medicine, Chicago, IL, USA \\
$^\ast$Corresponding author: \EMAIL{alikaan@u.northwestern.edu}}
}

\ABSTRACT{Unmanned aerial vehicles, commonly known as drones, have emerged as a disruptive technology with the potential to revolutionize operations across various industries. Drones are the fast-growing internet-of-things technology and are estimated to have a \$100 billion market value in the next decade. Exploring drone operations through research has the potential to yield innovative academic insights and create significant practical effects in diverse industries, offering a competitive edge. Drawing insights from both academic and industry literature, this article describes how technological advancements in UAVs may disrupt traditional operational practices in different industries (e.g., commercial last-mile delivery, commercial pickup and delivery, telecommunication, insurance, healthcare, humanitarian, environmental, urban planning, homeland security), identifies the value of this evolving disruptive technology from sustainability and operational innovation perspectives, argues the significance of this area for operations management by conceptualizing a research agenda. The current state of the art focuses on the computing aspect of analytical models to tackle a variety of synthetic drone-related problems, with mixed integer optimization being the primary tool. There is a very significant research gap that should focus on drone operations management with industry know-how by partnering with actual stakeholders and using a variety of tools (i.e., econometrics, field experiments, game theory, optimal control, utility functions). This article aims to promote research on UAVs from operations management and industry-specific point of view.
}

\HISTORY{This paper was first submitted in December 2023.}

\maketitle

\vspace{-0.9cm}

``\textit{Drones overall will be more impactful than I think people recognize, in positive ways to help society.}" -Bill Gates

\vspace{-0.3cm}

\section{Introduction}

In recent years, drones and unmanned aerial vehicles (UAVs) have emerged as a game-changing technology with enormous potential in a variety of applications. From their initial development for military use, with numerous mechanical advancements being made to improve their capabilities and performance, the adoption of drones has transformed the way organizations from different industries operate, facilitating more efficient, cost-effective, and sustainable processes. The introduction of Beyond Visual Line of Sight (BVLOS) operations has enabled drones to cover larger distances and operate in inaccessible locations, a once heavily restricted practice due to safety concerns, now made possible through drone technology and regulatory progression. Alongside, autonomous flight capabilities have evolved, allowing drones to operate without human intervention and handle functions like takeoff, landing, and navigation. With technology improvement, the payload capacity of drones has also escalated, accommodating heavier and more complex payloads, which in turn has broadened applications. In addition, enhanced battery life has allowed longer flights over greater distances, overcoming past limitations. On the regulatory front, entities like the Federal Aviation Administration (FAA) have simplified the drone licensing process, ensuring legal and safe operations \citep{lieb2020unmanned}. For example, in 2019, the FAA granted approval to UPS Flight Forward, making it the pioneering drone service to operate as a commercial airline \citep{insider2020amazon}. Additionally, in January 2021, the FAA for the first time gave permission to American Robotics to operate fully automated drones without the presence of humans and noted that they are working on an Unmanned Aircraft Traffic Management System \citep{higgins2021drone}. In order to develop a regulatory baseline for drones in commercial and medical delivery, the UK allocated £90 million \citep{ukdot2020}. Moreover, the Indian government is actively investing in drone-related activities to make the country a global drone hub by 2030 \citep{kumas2022daas}. According to \citet{cornelldrone2023}, the direct operational expense for a drone delivery of a single package is approximately \$13.50, where the major cost is labor cost (i.e., piloting, training). With a diminishing need for extensive human control, the cost-benefit of using UAVs is expected to increase in the near future. Currently, in terms of carbon emission cost, UAVs are much more beneficial than electric cars, electric vans, and internal combustion engine vans, moreover, their increased autonomy will make drones much more cost-competitive. The continual technological advancements in UAVs will persist in altering the conventional operational practices of companies, giving rise to new business models and enhancing overall sustainability.

During the initial months of 2022, the global daily number of commercial drone deliveries surpassed 2,000 \citep{cornelldrone2023}. While the quantity of delivery drone units remains modest, with approximately 35,000 units as of the end of 2022, the growth rate is 55\% \citep{insider2020amazon}. Furthermore, companies such as Amazon (i.e. Prime Air), Walmart, Google, FedEx, UPS, DHL, Walgreens, Zipline, and Domino's, are actively investing in drone delivery services, with the potential to revolutionize the logistics industry \citep{macias2020optimal}. Walmart is emerging as a front-runner in commercial drone delivery, having completed more than 6,000 drone deliveries in 2022 \citep{insider2020amazon}. The drone delivery market in the US is projected to witness substantial growth, expanding from 1 billion US dollars to 9 billion US dollars between 2022 and 2030, encompassing platforms, software, infrastructure, and services \citep{statistadrones}. In addition, drones in the transportation sector alone have the potential to curtail carbon emissions by up to 4.5 billion tons annually \citep{agenda2016shaping}, which is a crucial step to fulfilling the United Nations Sustainable Development Goals. Aside from the commercial sector, drones are being extensively utilized in other sectors like healthcare and humanitarian operations. For instance, drones are being tested for organ transplantation \citep{sage2022testing} and they have been used for disaster response activities like aid delivery and damage assessment, since 2005 \citep{chowdhury2017drones}.

\begin{figure}
\begin{center}
\caption{Distribution of drone operations articles according to sectors and years.} \medskip
\includegraphics[width=0.75\textwidth]{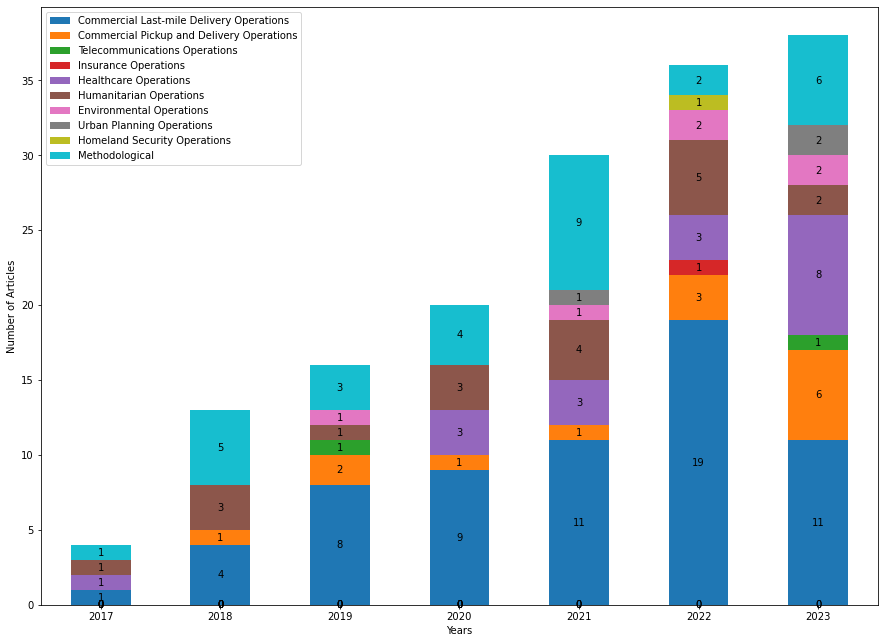}
\label{Img:Years}
\end{center}
\end{figure}

\begin{figure}
\begin{center}
\caption{Keyword analysis of drone operations articles ($n=158$).} \medskip
\includegraphics[width=0.9\textwidth]{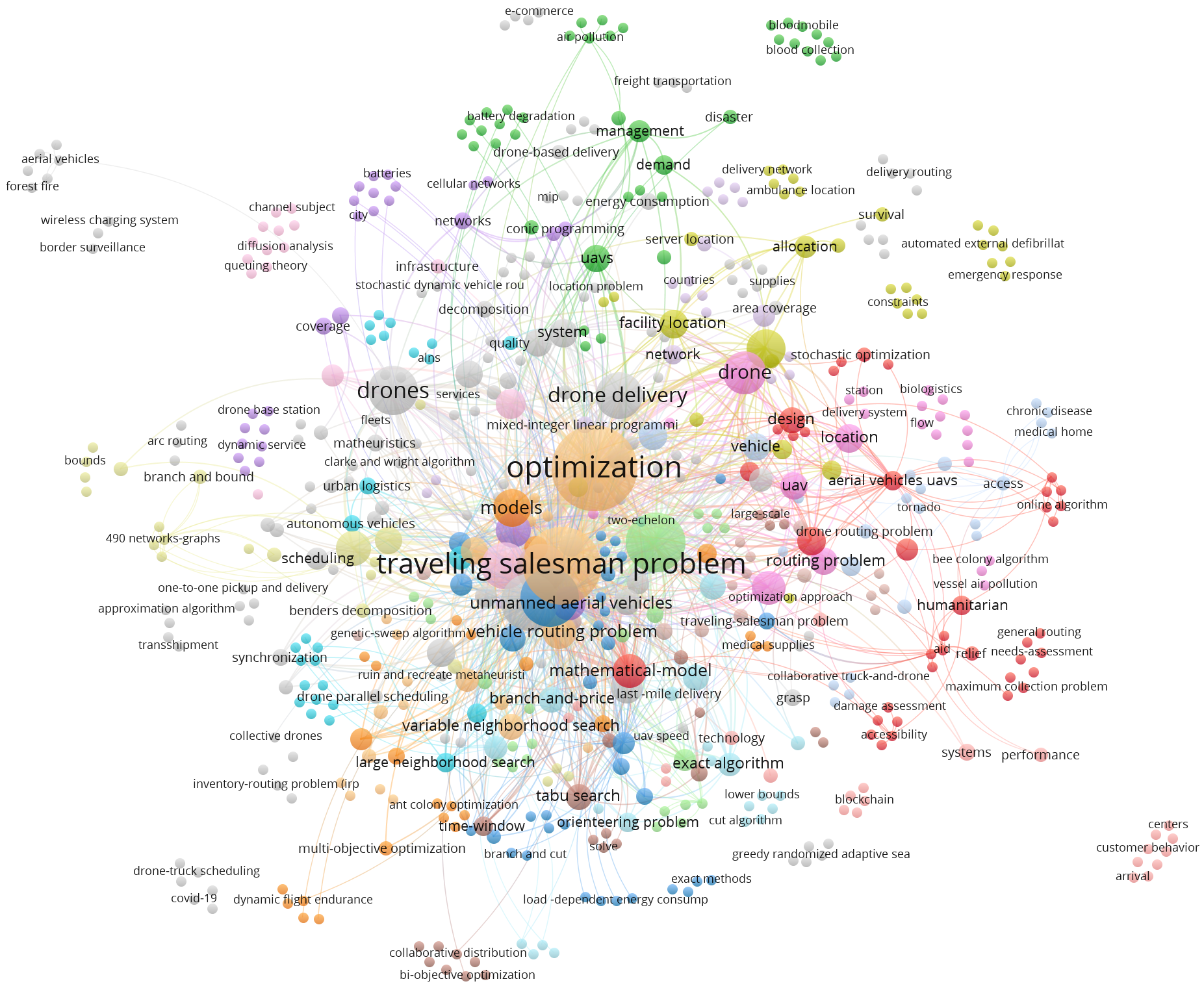}
\label{Img:Keywords}
\end{center}
\end{figure}

This paper aims to highlight the disruptive potential of unmanned aerial vehicles across various industries, emphasizing the need for research and industry collaboration to explore the operational and sustainability aspects of this technology. Our search methodology involved identifying relevant articles published or available online between 2017 and 2023 that contained an analytical model (i.e., optimization, queuing, Markov decision process). Through this approach, we identified a total of 158 articles that met the inclusion criteria: 64 commercial last-mile delivery, 14 commercial pickup and delivery, 2 telecommunications, 1 insurance, 18 healthcare, 19 humanitarian, 6 environmental, 3 urban planning, 1 homeland security, and 30 methodological. The articles in the methodological class develop solution algorithms for the problems that were already introduced. Distribution of drone operations articles according to sectors and years are depicted in {Figure \ref{Img:Years}}. The keyword analysis of articles, using \textit{VOSviewer}, is presented in Figure \ref{Img:Keywords}. The journal distribution (Table \ref{table:journals}) and keyword analysis for each industry (Figure \ref{Img:last_mile_keywords} - \ref{Img:other_keywords}) are presented in Appendix \ref{sec:appendix}.

To the best of our knowledge, this is the first article that examines drone operations from an industry, sustainability (i.e., social, environmental, economic), and operational innovation framework. We highlight the importance and potential of drones in OM/OR/MS research. We conclude that the research often centers on computational aspects, with a focus on large instances and synthetic data. Real-life cases are less frequent, and equity measures in drone operations are lacking. Furthermore, areas such as OM interfaces, procurement processes, and effects on the workforce remain unexplored. Some industries, like agriculture, construction, and green energy, lack analytical attention to drone applications. Integrating industry knowledge for practical applications is crucial. Lastly, the potential impacts of drones on supply chain financing, risk, and resilience warrant further investigation.

The remainder of this paper is organized as follows. The analytical model applications of drone technology in various industries are presented in {Section \S \ref{sec:Sectors}}. {Section \S \ref{sec:Sustainability}} includes the impact of drone operations on the United Nations Sustainable Development Goals. {Section \S \ref{sec:Innovation}} describes how the technological advancements in drones are leading to new innovative business models. Finally, we discuss several research directions in {Section \S \ref{sec:Conclusion}}.

\section{Drone Technology in Different Industries} \label{sec:Sectors}

In this section, our initial step involves the creation of a comprehensive framework for drone operations. This framework delineates key elements such as the primary hierarchy, network configuration, goods transported, synchronization methods, drone mobility, as well as the number of vehicles and drones. Subsequently, in the subsequent section, we utilize this framework to scrutinize the utilization of drones across diverse industries. Our analysis aims to incorporate sustainability practices and operational advancements within these applications.

\subsection{Echelons and Network Types}

\begin{figure}
\begin{center}
\caption{Drone operations framework.} \medskip
\includegraphics[width=0.9\textwidth]{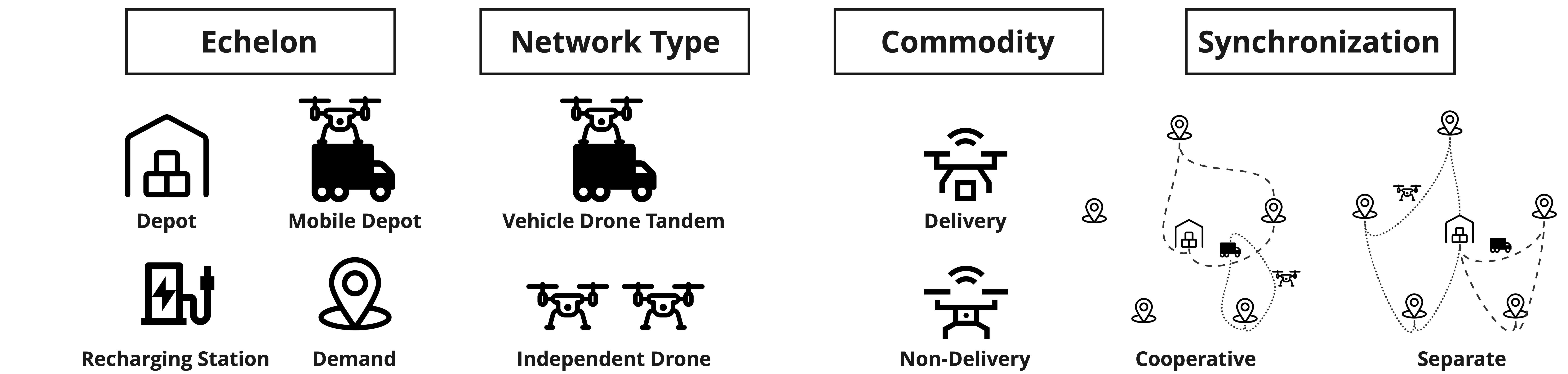}
\label{Img:drone-framework}
\end{center}
\end{figure}

In the context of drone operations modeling, an echelon framework is used across different settings. We identify four main echelons as: (i) \textit{supply node}, (ii) \textit{mobile platform}, (iii) \textit{recharging station} and (iv) \textit{demand node}. In a commercial context, the \textit{supply node} often takes the form of a warehouse or a distribution center, functioning as the primary hub where vehicles collect units slated for delivery. Meanwhile, in non-commercial settings, this depot could be a hospital or a relief base station. From these depots, fully-charged drones equipped with packages are dispatched to fulfill delivery needs. In cases where the drone is used for assessment, the \textit{supply node} is just a drone base. On the other hand, \textit{mobile platform} is crucial in vehicle-drone tandem systems. The drone remains with this \textit{mobile platform} until it is deployed, laden with a package, to cater to a demand node. Post-delivery, the drone can return to either a stationary depot or the \textit{mobile platform} to recharge its battery and/or collect more packages for delivery. In certain instances, the \textit{mobile platform} itself can also act as a service point, delivering packages directly to customers. In non-delivery operations, the \textit{mobile platform} is used to extend the coverage of the drones. \textit{Recharging stations} are specialized facilities designed solely to enable drones to recharge or swap their batteries, thereby extending their service range. The final echelon is the \textit{demand node}. This entity can take various forms depending on the context. It could be an area requiring surveillance or assessment, a specific location for meal delivery, or a more generalized point where drones deliver their packages. This echelon represents the endpoint of the delivery chain, effectively completing the cycle of drone-based operations.

Based on the echelon structure, drone network types predominantly concentrated on two main types of systems: \textit{vehicle-drone tandem systems} and \textit{independent drone systems}. \textit{Vehicle-drone tandem systems} are predicated on the cooperative interaction between a bigger vehicle, with higher endurance but less flexibility, and a drone. This synergistic system capitalizes on the respective strengths of bigger vehicles and drones, resulting in operations that are characterized by increased efficiency and reduced costs. For instance, while drones may be used for rapid delivery over difficult terrains, land and water vehicles may be employed to cater to high-volume or long-distance deliveries. On the other hand, \textit{independent drone systems} involve the exclusive use of drones for operations. These systems typically require comprehensive models to address a myriad of logistic considerations, such as drone flight scheduling, battery management, and load optimization, among others. Collectively, these two types of systems encapsulate the broad spectrum of operational frameworks that are central to the current state of drone operations.

There are certain actions that a drone can perform. In the context of \textit{delivery}, they can deliver or pickup goods, spray or sprinkle substances in the form of liquid or powder, and transmit signals. In contrast in the context of \textit{nondelivery}, with the camera or sensors attached to them, they can assess, inspect, and monitor. Finally, to sustain their batteries they can swap batteries or get charged in docs. Also, some drones can get closer to recharging stations to get charged in a contactless way. There are many factors affecting the energy usage of a drone such as its speed, weight of the payload, and wind.

In the context of \textit{vehicle-drone tandem systems}, synchronization refers to the ability of vehicle(s) and drone(s) to work together in a coordinated and synchronized manner to achieve a common goal. This phenomenon is the fastest growing internet-of-things (IoT) technology. The vehicle boasts a significant carrying capacity but typically operates at lower speeds in urban settings, while the drone, although faster and free from street network restrictions, has limited range and carrying capabilities. \citet{moshref2021comparative} show that a higher level of synchronization reduces lead times. Synchronization can be classified in terms of vehicle type, synchronization type, number of vehicles, number of drones, and nodes per drone trip. The vehicles (\textit{mobile platform}) can be truck (88 articles), mothership drone (6 articles), ambulance (3 articles), ship (2 articles), bus (1 article), and motorcycle (1 article). In \textit{cooperative synchronization}, the drone takes off and lands on the vehicle. In contrast, in \textit{separate synchronization} the vehicle and the drone perform the task separately to increase coverage. In cases where the drone visits both the vehicle and a stationary point (i.e. depot, en-route), we classify these synchronizations as both. In the fleets, there can be single or multiple drones and vehicles. The drones can visit single or multiple nodes after taking off from a vehicle. Finally, in some cases, the drone can/must loop, where it visits and comes back to the vehicle, while the vehicle is stationary. For more information on detailed features of each drone operations article refer to {Table \ref{table:detailed-features}} in {Appendix \S \ref{sec:appendix}}.

The primary types of decisions addressed are routing (129 articles), allocation (73 articles), location (43 articles), distribution (39 articles) and inventory (29 articles). Routing decision (operational) covers decisions on vehicle/drone movements and visiting prioritization. Allocation (operational) refers to allocating drones to stations, demand to drones, or scheduling multiple drones in a single vehicle by allocating a drone to time slots. Location (strategic) encapsulates deciding the optimal position of a facility, such as \textit{recharging station}. Distribution (operational) is about decisions regarding the distribution of multiple goods to demand nodes. Finally, inventory decisions (tactical) manage various aspects of the inventory. Note that for drones with unit capacity, distribution and inventory decisions are not considered, since it is handled through routing variables. In such cases, when a drone visits a demand node, it is assumed that the unit good is delivered. For more information on decision levels of each drone operations article, along with objective functions, refer to {Table \ref{table:DL-obj}} in {Appendix \S \ref{sec:appendix}}.

\subsection{Applications in Different Industrial Sectors} \label{subsection:sectors}

In this section, we categorize and examine drone operations across various sectors within commercial last-mile delivery (Section \S \ref{subsec:221}), commercial pickup and delivery (Section \S \ref{subsec:222}), telecommunication (Section \S \ref{subsec:223}), insurance (Section \S \ref{subsec:224}), healthcare (Section \S \ref{subsec:225}), humanitarian efforts (Section \S \ref{subsec:226}), environmental applications (Section \S \ref{subsec:227}), urban planning (Section \S \ref{subsec:228}), and homeland security (Section \S \ref{subsec:229}). Please note that our focus here is solely on analytical studies. While there may be additional application domains, those lacking analytical modeling are excluded from this discussion, and remain as future research directions. 

\subsubsection{Commercial Last-mile Delivery Operations} \label{subsec:221}

Drones are mainly used in last-mile delivery in e-commerce. With increased competition in the e-commerce sector, companies are using drones to sustain their on-demand and same-day delivery operations. Due to the COVID-19 pandemic, the volume of e-commerce proliferated and people started living in dispersed areas. To maintain customer satisfaction in all geographical regions, lead times should be at a minimal level. A substantial body of literature in the field of commercial last-mile delivery under \textit{vehicle-drone tandem systems} is built upon the seminal work of \citet{murray2015flying}. This work introduces the flying sidekick traveling salesman problem (FSTSP), where a truck embarks on its route from a designated depot, carrying a load of parcels destined for delivery at various nodes. Alongside these parcels, the truck is equipped with a drone, which is also capable of delivering a single package to a demand node at a time. This introduces a level of flexibility, efficiency, and sustainability previously unattainable with traditional delivery systems. Models that tackle the FSTSP are charged with the task of addressing a multitude of interrelated elements. Among these elements are the routing and scheduling of the truck, the identification of optimal takeoff and landing points for the drone, and the effective management of the drone itself. Note that all articles in this sector assume there is a single depot.

There are articles that extend the single vehicle, single drone system, which was introduced in \citet{murray2015flying}. \citet{campuzano2023drone} consider different flight times for the drone depending on the selected speed, payload weight, and weather conditions. All of these factors also affect the energy consumption of the drone. Accounting for the energy expenditure of drones in the model ensures that the flight plans generated are sustainable and would not lead to premature drone grounding due to battery drainage, thus significantly improving the overall system reliability. Moreover, \citet{dell2021drone} allows drones to wait in the customer node for improved synchronization. \citet{jeong2019truck} consider the effect of parcel weight on drone energy consumption and introduce restricted flying areas. The inclusion of such constraints ensures the respect of regulatory boundaries and further enhances the safety and legality of drone operations. \citet{liu2022flying} incorporate stochastic travel times. Additionally, in \citet{masone2022multivisit}, the drone is capable of carrying multiple packages, thus it visits multiple demand nodes when launched from the vehicle. Finally, \citet{najy2023collaborative} embed inventory routing problem into FSTSP.

Another technological advancement includes the utilization of multiple drones in a single vehicle, which enhances the efficiency of the system through parallel processing capabilities \citep{murray2020multiple}. There are articles that extend the single vehicle, multiple drone system. \citet{boysen2018drone} schedule the delivery to customers by drones for given truck routes. \citet{bruni2022logic} allows drones to do multiple trips while the vehicle is stationary at a node. \citet{chang2018optimal} focus on vehicle routing decisions to make wider drone-delivery areas using machine learning. In addition, \citet{kang2021exact} consider a heterogeneous drone fleet, where the drones have different speeds and battery capacities. In \citet{luo2021multi} and \citet{poikonen2020multi}, drones are capable of carrying multiple packages, thus they visit multiple demand nodes when launched from the vehicle. Since in big cities people do not prefer to go to the ground floor to receive packages, \citet{momeni2023new} consider drone delivery to different altitudes. Moreover, \citet{raj2020multiple} consider variable drone speeds. The speeds are inversely proportional to drone energy consumption. Although the majority of the works consider a truck as a vehicle, \citet{wang2022piggyback} is the only article in this class to consider a mothership drone. This work assumes that the dispatched drones go back to the stationary depot after delivery. \citet{yang2023planning} include traffic uncertainty assuming that as drones are employed more in different sectors, it will result in air traffic. Finally, \citet{zhao2022robust} include navigation environment uncertainty and synchronization risk. 

Another substantial portion of the articles consider multiple vehicles and multiple drones. This problem is known to be the vehicle routing problem with drones (VRP-D) in the literature. Despite, in the majority of the \textit{cooperative synchronization} works every vehicle has its own dedicated drone fleet, \citet{bakir2020optimizing} allow a drone dispatched from one truck to land on another truck. \citet{chen2022deep} consider \textit{separate synchronization} under demand uncertainty for same-day delivery operations using deep learning. \citet{coindreau2021parcel} include time-window constraints for modeling on-demand delivery. Furthermore, \citet{kitjacharoenchai2019multiple} consider both synchronization types where on top of a vehicle-drone tandem system, drones can directly travel from the depot to customers or drones that are dispatched from vehicles can fly back to the depot. \citet{kitjacharoenchai2020two} and \citet{masmoudi2022vehicle} extend the VRP-D, drones have the ability to transport multiple packages, allowing them to visit multiple demand nodes upon being launched from the vehicle. Additionally, \citet{kloster2023multiple} consider the location decision of robot-supported packet stations, where drones can grab a package. \citet{nguyen2022min} and \citet{saleu2022parallel} consider \textit{separate synchronization} of trucks and drones, where the drones can loop from depot to a customer. In these articles, drone scheduling and vehicle routes are optimized. \citet{rave2023drone} include deciding on locations for dedicated drone stations, called microdepots, from where drones can start as well as deciding on the fleet of conventional trucks, trucks equipped with drones, and drones. \citet{ulmer2018same} consider same-day delivery with a heterogeneous drone fleet and demand uncertainty. In order to determine the delivery method for an order, they propose a policy function approximation based on geographical districting. This approach decides whether the order should be fulfilled using a drone or a vehicle. Their conclusion indicates that the implementation of geographical districting is highly effective in increasing the expected number of same-day deliveries. Furthermore, \citet{wang2022truck} considers time-dependent road travel time, due to traffic and proposes a model for improved synchronization. Finally, in \citet{zhang2022novel} the maximum flight endurance of drones is dynamically adjusted based on their loading rate to meet the requirements of practical application scenarios.

In order to increase drone flexibility, several articles consider en-route drone operations, where the drones can land, wait and relaunch from discrete locations on the arc that are non-customer nodes \citep{salama2022collaborative,schermer2019hybrid,thomas2023collaborative} or intermediary docking hubs \citep{wang2019vehicle,wen2022heterogeneous,xia2023branch}. This aspect allows for the dynamic selection of launch and retrieval points, thereby enhancing the adaptability and ability of the system to respond to unpredictable changes in the operational environment. Another concept in commercial last-mile delivery operations under \textit{vehicle-drone tandem systems} for improved on-demand and same-day delivery is drone resupply \citep{dayarian2020same,dienstknecht2022traveling}. In drone resupply, the truck is routing customer nodes and when the inventory of trucks decreases, the drone resupplies the truck from the depot. Using drones for resupply can reduce the total delivery time by up to 20\% \citep{pina2021traveling}. Despite the fact that many of the articles in this class aim to minimize cost or minimize the makespan, there are some objectives that minimize carbon emission or energy usage for environmental sustainability \citep{chiang2019impact,coindreau2021parcel,dukkanci2021minimizing,zhang2022novel}.

Modeling research has also considered \textit{independent drone system} for commercial last-mile delivery operations. In these systems, to support coverage of drone delivery, models decide where to strategically locate intermediate \textit{recharging stations} \citep{cokyasar2021designing,huang2020method,pinto2022point,shavarani2019congested} and drone dispatch centers. Government regulations, customer behaviors \citep{baloch2020strategic}, gridding choices \citep{betti2023dispatching}, energy functions \citep{cheng2020drone}, and coverage types \citep{liu2023elliptical} influence these location decisions. These location decisions play a pivotal role in determining the overall efficiency and effectiveness of the drone network, and their selection constitutes a significant part of network design considerations, especially recharging stations allow for the continual operation of drones and provide a much-needed buffer against energy depletion, thus contributing to the development of a more resilient and robust drone network. Modeling research has further delved into the issue of battery allocation across these carefully selected recharging stations, exploring ways to optimize the distribution of power resources within the network \citep{asadi2022monotone}. Additionally, the uncertainties modeled in this class are demand \citep{asadi2022monotone,chen2021improved,cokyasar2021designing,shen2021operating,liu2022energy}, waiting time \citep{shavarani2019congested}, and wind direction \citep{jung2022drone}. Also \citet{wang2023optimal} consider a heterogeneous drone fleet. Although the majority of the works either try to minimize cost or distance, there are some works that use environmental sustainability objective functions, where the models minimize energy consumption \citep{cheng2020drone,liu2022energy,troudi2018sizing,xia2021joint}. To operate in an environmentally sustainable manner, in \citet{liu2022energy} drones have adjustable flight speed to trade off energy efficiency.

Innovative models include the application of blockchain technology for the management of shared drone fleets \citep{xia2021joint}. Blockchain presents promising opportunities for enhanced coordination and utilization of shared drone resources. Another scholarly work introduces the innovative notion of a collective drone system \citep{nguyen2023parallel}. This unique approach allows multiple drones to collaboratively participate in the delivery of a package, particularly when the weight of the parcel exceeds the carrying capacity of any individual drone. This concept showcases a pivotal enhancement in tackling weight-based limitations, thereby extending the operational capabilities of drone delivery systems. \citet{shen2021operating} consider a multi-depot drone delivery system where they test different allocation scenarios. They conclude that the shared allocation of drones rule gives a larger throughput capacity than the dedicated drone allocation to warehouses rule, and it reduces the operating cost. As drone technology advances and becomes more cost-efficient, the conclusion drawn is that independent drone delivery networks will experience greater decentralization, facilitating quicker product deliveries \citep{perera2020retail}.

\subsubsection{Commercial Pickup and Delivery Operations} \label{subsec:222}

Commercial pickup and delivery operations involve the use of drones to transport goods from one location to another, providing fast and efficient delivery services. Drones travel from drone hubs to pick up packages from a designated starting point and deliver them to specified destinations, streamlining the delivery process and reducing human intervention in the logistics chain. Unlike \textit{commercial last-mile delivery operations}, in pickup and delivery, there are multiple supply nodes that are dispersed. However, the aim of the models is the same: to sustain high satisfaction levels for on-demand and same-day deliveries, to a dispersed set of customers. In \textit{commercial pickup and delivery operations} articles there are models that consider \textit{vehicle-drone tandem systems} \citep{cheng2023adaptive,gu2023dynamic,ham2018integrated,jiang2023multi,karak2019hybrid,luo2022last,meng2023multi} and \textit{independent drone systems} \citep{gomez2022pickup,kyriakakis2023grasp,levin2022branch,liu2019optimization,liu2023routing,pei2021managing,pinto2020network}. Regardless of network type, many of the articles consider time-window constraints to model on-demand nature of the problem and energy aspects (i.e. battery swaps, weight, and distance based energy consumption) of the drones. Additionally, \citet{gu2023dynamic} considers uncertain demand arrival rate and dynamic vehicle routing. \citet{jiang2023multi} consider flexible docking where drones are not assigned to specific trucks and land on another truck. They also consider simultaneous pickup and delivery, for cost savings. Furthermore, \citet{levin2022branch} consider urban airspace, where the sky over current roadways functions as an airspace, structured akin to a multi-tiered transportation network. Finally, a novel work by \citet{cheng2023adaptive} suggests combining passenger and goods transportation through on-demand buses and drones, offering an innovative approach compared to traditional methods that rely solely on ground vehicles. This proposal involves utilizing on-demand buses for both passenger and parcel delivery, complemented by drones exclusively designed for parcel transportation, showcasing a unique integration strategy. This is the only study, where ground transportation is done through public transportation vehicles, with designated routes. The model aims to optimize both passenger and package delivery by using public transportation vehicles as a \textit{mobile platform}.

While the majority of research in this industry predominantly discusses drone utilization within the broad context of package or parcel delivery, some of \textit{independent drone system} articles tailor their studies towards meal delivery. These studies consider the operational dynamics unique to this industry, including factors such as clustering meals as being hot or cold for serving on-demand meal deliveries using dynamic routing \citep{liu2019optimization}, utilizing uncertain demand arrival to optimize price and fleet size of drones shared by multiple restaurants \citep{pei2021managing}, and including customer wait times, the geographical distribution of the population, the variety of restaurants available for network design \citep{pinto2020network}. The objectives considered in these articles are maximizing coverage or minimizing distance between customers and restaurants. In order to expand the coverage \citet{liu2019optimization} and \citet{pinto2020network} consider locating \textit{recharging stations}.

\subsubsection{Telecommunication Operations} \label{subsec:223}

In the context of telecommunication operations, drones are used for broadband internet delivery \citep{cicek2019location,colajanni2023three}. \citet{cicek2019location} focuses on the 3D location problem of multiple drone base stations (DBSs) and resource allocation in a wireless communication network. The problem is formulated as a dynamic capacitated single-source location-allocation problem, with the capacity of a DBS and data rate being non-linear functions of distance and allocated resources, resulting in a non-linear model. On the other hand, \citet{colajanni2023three} introduces a three-stage stochastic network-based optimization model, under demand uncertainty, for providing 5G services using UAVs in disaster management phases. It involves ground users/devices requesting services from controller UAVs, which are executed by a fleet of interconnected UAVs functioning as a flying ad-hoc network via 5G technology.

\subsubsection{Insurance Operations} \label{subsec:224}
Drones are utilized in insurance operations for risk assessment, property monitoring, and claims inspection, providing insurers with quick and accurate data to streamline processes, improve customer service, and enhance risk management. As an analytical model, \citet{zeng2022nested} discusses the routing and coordination of a drone-truck pairing for specified observation tasks at multiple locations, with periodic battery swaps, and presents a case study where an insurance company used this system to inspect and collect damage evidence after a fire disaster from 631 clients' properties, in the minimum amount of time.

\subsubsection{Healthcare Operations} \label{subsec:225}

In the context of healthcare operations, drones are used for medical good distribution \citep{enayati2023multimodal,gentili2022locating,enayati2023vaccine,ghelichi2021logistics,macias2020optimal,park2022scheduling,wang2023multi,zhang2023collaborative,wang2023robust}, and medical good collection \& distribution \citep{dhote2020designing,hou2021integrated,kim2017drone,rezaei2023integrating,shi2022bi}. Drones equipped with medical supplies and equipment can reach remote areas quickly, delivering life-saving vaccines \citep{dhote2020designing,enayati2023multimodal,wang2023robust,enayati2023vaccine,zhang2023collaborative}, blood \citep{dhote2020designing,hou2021integrated,rezaei2023integrating}, testing kits \citep{kim2017drone,park2022scheduling}, and medications \citep{kim2017drone} to those in need. The ability of drones to quickly and efficiently deliver these supplies has been particularly important in remote areas during emergencies. For example, during the COVID-19 pandemic, drones were used to deliver test kits, masks, and other medical supplies to areas that were difficult to reach by ground transportation. Drones can also pickup and deliver medical samples to laboratories for testing, allowing for faster diagnosis and treatment. Many medical supplies may have strict storage and handling requirements to ensure their viability. This necessitates the incorporation of perishability into models \citep{enayati2023multimodal,gentili2022locating,zhang2023collaborative}. The most common uncertainty types in this context is demand \citep{hou2021integrated,rezaei2023integrating,wang2023multi}, supply \citep{rezaei2023integrating,wang2023robust} and drone endurance \citep{wang2023multi}. The majority of the analytical models employ a multi-objective approach to minimize operational costs and ensure equity or social sustainability. For example \citet{macias2020optimal} minimizes excess vaccine delivery among locations to ensure geographical equity similarly, \citet{wang2023multi} minimizes maximum response time. For social sustainability, \citet{wang2023robust} maximizes the number of successful vaccinations.

Drones are increasingly being used for emergency response in healthcare, particularly in situations where time is critical. For example, the use of drones has been explored to deliver automated external defibrillators (AEDs) to people in cardiac arrest \citep{boutilier2022drone,lejeune2022maximizing,shin2022heterogeneous,wankmuller2020optimal}, reducing response times and increasing survival rates. For a bystander's use of an AED, \citet{boutilier2022drone} and \citet{shin2022heterogeneous} synchronize the drone carrying the AED with the established ambulance system. The goal here is not to replace the traditional ambulance system but to complement it, ensuring that AEDs can be delivered rapidly to the scene of a cardiac emergency while the ambulance is en route. This strategy optimizes the use of resources, potentially shortening the critical window between cardiac arrest and defibrillation, thus increasing survival rates. It is important to understand that while drone-to-ambulance synchronization represents a significant technological advancement, it also requires the careful coordination of human elements (bystanders) to ensure the optimal functioning of this innovative service. The most common uncertainty types in this context are demand arrival rate \citep{lejeune2022maximizing,shin2022heterogeneous} and response time \citep{boutilier2022drone}.

\subsubsection{Humanitarian Operations} \label{subsec:226}

In the context of humanitarian operations, drones are used for relief distribution \citep{chowdhury2017drones,dukkanci2023drones,ghelichi2022drone,jeong2020humanitarian,kim2019stochastic,rabta2018drone,yin2023robust,zhu2022two}, and relief collection \& distribution \citep{lu2022multi}. The unique characteristics and requirements inherent to the specific disaster scenarios of these analytical models often necessitate the consideration of additional factors, compared to their commercial counterparts. For instance, \citet{dukkanci2023drones} and \citet{yin2023robust} consider earthquake, \citet{chowdhury2017drones} and \citet{ghelichi2022drone} consider hurricane, and \citet{jeong2020humanitarian} consider a civil war. Given the extensive variety of aid products required in disaster relief, these models generally do not specify a particular group of items. Also, \citet{rabta2018drone} uses demand node prioritization for scarce resource management. \citet{jeong2020humanitarian} introduce the concept of humanitarian flying warehouse (HFW), where the HFW is an airship stationed at high altitudes and acts as a mothership drone. The authors propose that this system, with its potential to efficiently traverse geographical barriers and expedite aid distribution, may well be the catalyst for a new era in humanitarian aid logistics. Note that all of the problem settings can be classified as \textit{response operations} in terms of the well-known \textit{disaster management cycle} framework. The most common uncertainties in humanitarian delivery problems are demand \citep{dukkanci2023drones,ghelichi2022drone,yin2023robust,zhu2022two}, drone flight distance \citep{kim2019stochastic} and drone flight time \citep{yin2023robust}.

In addition to the humanitarian \textit{delivery} models outlined, \textit{non-delivery} models primarily involve assessment, surveying, and identification activities. Drones in these contexts may be used for tasks such as post-disaster damage assessment, mapping, and surveying disaster-affected areas, identifying individuals or groups in need of search \& rescue, or other similar tasks where the value of the drone lies in its ability to gather and relay information, in a flexible way. Elaborating further, the drones can be used for network assessment, where they are deployed to evaluate a particular road or logistics network \citep{escribano2020endogenous,reyes2021exploration,zhang2021humanitarian,zhang2023robust}. This evaluation is occasionally supplemented by the use of ground vehicles to ensure comprehensive and accurate assessment \citep{oruc2018post,zhang2021humanitarian}. These models are crucial in situations where a network's operability and integrity are compromised due to disasters, and there is a need to understand the extent of damage for restoration efforts. Conversely, drones can be used in area assessment which, involves the use of drones to survey disaster-stricken regions \citep{bravo2019use,chowdhury2021drone,glock2020mission,grogan2021using,kyriakakis2022cumulative,ozkan2023uav}. The main objectives of such models can range from damage assessment of the region \citep{chowdhury2021drone,grogan2021using} to aiding search \& rescue operations \citep{grogan2021using,kyriakakis2022cumulative,bravo2019use}. The versatility of drones allows them to access regions otherwise difficult or hazardous for humans, thereby increasing the safety and efficiency of these operations. A notable consideration in one study involves the integration of priority levels for specific areas \citep{bravo2019use}. This factor aids in the strategic allocation and efficient use of drone operation resources, ensuring that high-priority areas receive attention commensurate with their immediate needs. The most common uncertainties in humanitarian non-delivery problems are road network \citep{escribano2020endogenous,reyes2021exploration} and response time \citep{zhang2023robust}.

\subsubsection{Environmental Operations} \label{subsec:227}

Drones are instrumental in environmental monitoring and conservation efforts. Equipped with sensors, they enable efficient data collection, aiding in wildlife tracking, biodiversity assessments, and ecosystem management. Drones provide valuable insights into deforestation, pollution levels, and climate change impacts, facilitating evidence-based decision-making for sustainable practices. In the context of analytical models, drones are used for air pollution detection. According to the International Maritime Organization (IMO), shipping contributes to approximately 13\% of sulfur dioxide, nearly 15\% of nitrogen oxides, and almost 2.7\% of carbon dioxide emissions from human activities \citep{lindstad2021reduction}. In response to IMO regulations, drones are now widely employed for monitoring ship emissions in ports \citep{liu2023exact,shen2020synergistic,shen2022synergistic,xia2019drone}. While in \citet{liu2023exact} and \citet{shen2022synergistic} drones work in tandem with ships, \citet{shen2020synergistic} and \citet{xia2019drone} consider independent drone systems. Note that these works consider moving targets. The ships that are being monitored are moving in the ports. To include this aspect of the problem \citet{liu2023exact} consider heterogeneous drone fleets with different flight speeds and endurance and \citet{xia2019drone} prioritize inspection of vessels with higher weight.

Drones can also be utilized for forest fire risk mitigation and precision agriculture. \citet{ozkan2023uav} introduce simulation-based optimization models to create routes for drones that detect forest fires, particularly in remote areas distant from residential regions. The models account for uncertainties in forest fire occurrences and drone flight speeds. Moreover, precision agricultural systems offer data-driven and scientific assistance for managing crops. \citet{ermaugan2022learning} presents a learning-based algorithm capable of solving the drone routing problem in real time, with the flexibility to update routing decisions whenever operational conditions change.

\subsubsection{Urban Planning Operations} \label{subsec:228}

The most common urban planning drone operation is traffic inspection. \citet{amorosi2023multiple} address optimization problems related to coordinating a tandem of a helicopter and a fleet of drones. The drones are launched from the mothership to perform specific tasks, with a particular focus on preventing and identifying potential concentrations of people during events like popular or religious festivals, taking into account the constraints imposed by COVID-19 restrictions. Similarly, \citet{campbell2021solving} introduces a continuous optimization problem, where a fleet of homogeneous drones needs to traverse a network comprising curved or straight lines. The drones have flexibility, as they can enter a line from any point, service a portion of it, exit through another point, and directly travel to any point on another line. Aside from traffic monitoring, this problem is applicable to inspecting pipelines, railways, or power transmission lines. Additionally, drones can be utilized in urban disinfection. For instance, \citet{manshadian2023synchronized} investigate the use of spraying UAVs to combat environmental contamination and infections during crises, like epidemics, with the help of ground vehicles, equipped with disinfection solutions. 

\subsubsection{Homeland Security Operations} \label{subsec:229}

Manned systems for border monitoring are risky and costly due to the borders' vastness and harsh conditions. A potential solution is employing small drones for smart border patrol to access inaccessible areas, improve response time, and enhance agent safety. To address the limited flight duration caused by battery constraints, \citet{ahmadian2022smart} suggests a continuous border surveillance system using drones with a dynamic wireless battery charging system along an electrification line, with flight intervals managed based on border section criticality. Their analytical models aim to minimize the number of drones employed, as well as, the total length of the e-line system.

\subsection{Methodological Aspects}

The analytical modeling approaches, uncertainties and solution approaches for each drone operations article are presented in Table \ref{table:models-uncertainties-solutions} in Appendix \S \ref{sec:appendix}. Note that there is a \textbf{Methodological} articles class. These articles develop computational methods to solve problems that were already introduced by other articles and, hence, do not contribute to the application side of the literature. In terms of modeling approach, 145 articles employ optimization modeling (i.e., discrete (118), conic (8), stochastic (8), robust (7), dynamic (3), online(1)), 16 articles employ stochastic process modeling (i.e., queuing (9), MDP (6), simulation (1)) and 6 articles employ other modeling (i.e., real analysis (2), constraint programming (1), continuous approximation (1), fuzzy programming (1), utility functions (1)) approaches. On the other hand, to solve these analytical models, there are 55 articles that develop exact solution methods, in contrast, the remaining 103 articles employ heuristic solution methods. Among these articles 33 of them consider at least one uncertain parameter, and the remaining 125 articles consider full information. In order to test the quality of the solution methods, 106 articles use synthetic instances, whereas 52 use real-life case studies. Under deterministic parameters, the computational complexity of the problem increases exponentially as the number of nodes increases, and since routing decisions are combinatorially complex, the exact solutions perform poorly in large instances. They are generally able to solve up to 50 nodes. Data-driven heuristic models are more dominant in the literature to be able to solve large instances with uncertainty involved. In short, the literature highly focuses on the computing aspect.

\section{Drone Technology and United Nations Sustainable Development Goals} \label{sec:Sustainability}
In this section, we explore the impact of drone technology on United Nations Sustainable Development Goals (SDGs) \citep{kitonsa2018significance}, by highlighting their potential to accelerate progress and address key challenges through operations \citep{besiou2021humanitarian,sunar2022socially}. A full list of each of the 17 SDGs is provided in Table \ref{table:UNSDG} in Appendix \ref{sec:appendix}. We classify the SGDs under three categories: social sustainability (\textit{SDGs 1, 2, 3, 4, 5, 10, 11, 16, 17}), environmental sustainability (\textit{SDGs 6, 7, 13, 14, 15}) and economic sustainability (\textit{SDGs 8, 9, 12}).

\subsection{Drone Technology and Social Sustainability}

Drones have emerged as a transformative tool with the potential to significantly impact social sustainability through effective operations management. They facilitate last-mile delivery of essential goods such as food, water, medication, vaccine, and chemotherapy, to remote areas (\textit{SDGs 2, 3}). During the COVID-19 pandemic, drones became the primary tool for food delivery, transporting medical samples, and quarantine supplies in different countries. \citet{cozzens2020china} report that the utilization of drones has led to a transport speed enhancement of over 50\% in comparison to road transportation in China. This underscores their effectiveness as a more proficient mode of transport for the purposes of epidemic prevention and control. In addition, \citet{parkhill2022covid} report that during a rigorous COVID-19 lockdown in Shanghai, China, firefighters at the local level were employing drones to carry out contactless deliveries of essential items such as food and necessities. In terms of crowd post-lockdown monitoring, drones with thermal sensors can detect people with a fever in a crowd for public health \citep{mohsan2022role} (\textit{SDG 3}). They support sustainable agricultural productivity by real-time crop monitoring and enabling precision spraying (\textit{SDG 2}). They also promote economic opportunities in under-served communities, by increasing the accessibility of local businesses at a low cost (\textit{SDG 1}). Drones aid in urban planning and sustainable development by facilitating aerial mapping, supporting infrastructure maintenance, and improving disaster preparedness. They enable efficient urban planning, assist in monitoring environmental indicators, and enhance emergency response capabilities (\textit{SDG 11}). They also support peacekeeping and law enforcement by enhancing surveillance capabilities. They assist in border surveillance, provide situational awareness during emergencies, and support evidence collection for justice systems \citep{al2023systematic} (\textit{SDG 16}). These, in turn, promote social inclusivity and equitable resource provision, contributing to the overall betterment of society's well-being and resilience.

\subsection{Drone Technology and Environmental Sustainability}

UAVs offer a versatile range of applications that contribute to more efficient resource management, reduced ecological impact, and enhanced conservation efforts. Both in terms of carbon emission and energy consumption, drones are much more sustainable compared to many other vehicles \citep{goodchild2018delivery}. Drones support environmental sustainability in water, land, and air. Drones aid in water resource management by monitoring water quality, assessing water bodies, and identifying pollution sources \citep{sibanda2021application}. They enable efficient monitoring of water resources, support disaster response in water-stressed regions, and facilitate early detection of water-related issues \citep{restas2018water,zhai2023improving} (\textit{SDG 6, 14}). Drones contribute to the promotion of affordable clean energy by surveying potential sites for renewable energy installations. They are used for inspecting power lines, wind turbines, and solar farms, by high-definition cameras and thermal sensors capture detailed imagery, aiding in identifying maintenance needs \citep{addabbo2018uav} (\textit{SDG 7}). Drones aid in biodiversity conservation and land management by monitoring wildlife populations, combating poaching, and facilitating reforestation efforts. They support habitat mapping, aid in wildlife protection, and assist in monitoring land degradation. In agriculture, precision farming utilizing UAVs facilitates targeted pesticide and fertilizer application, minimizing waste and runoff while optimizing crop yields. In forestry, drones aid in the early detection of forest fires, enabling swift responses that mitigate widespread destruction. Additionally, UAVs are indispensable in wildlife monitoring, enabling researchers to gather essential data without disrupting natural habitats \citep{jimenez2019drones} (\textit{SDG 15}). On top of all of the mentioned benefits, drones contribute to climate action by providing data for climate research, monitoring deforestation, and assisting in disaster resilience planning. They support environmental monitoring, aid in the assessment of climate change impacts, and facilitate early warning systems for natural disasters (\textit{SDG 13}). In summation, UAVs constitute a multidimensional instrument for advancing environmental sustainability and nurturing a more resilient planet.

\subsection{Drone Technology and Economic Sustainability}

Drones foster sustainable economic growth and decent work opportunities by stimulating the drone industry and creating jobs. They enable the development of drone-related businesses and provide new employment avenues in various industries. A study by \citet{dronedeploy2020} states that the introduction of automation technologies will lead to the removal of 800,000 jobs requiring lower skills; however, it subsequently will generate 3.5 million fresh employment opportunities in the UK. These recently created positions will boast an average annual salary that surpasses the lost jobs by \$13,000 (\textit{SDG 8}). In addition, drones drive innovation by revolutionizing industries and support entrepreneurship through various industries, by allowing fresh business models. The disruptive technology enables both product and market growth \citep{giones2017toys} (\textit{SDG 9}). Finally, drones support responsible consumption and production by optimizing resource use. For instance, in agriculture, farmers are under pressure to discover sustainable methods of feeding the population as a result of the fast expansion of the global population. Farmers assess environmental factors that affect crop yield using information from drone technology. Drones are being utilized for precision agriculture which minimizes the wasted resource usage. Also, ranchers are able to monitor the health of their animals and implement efficiency measures that can lessen their impact on the environment, increase their production, and save costs \citep{pinguet2021role} (\textit{SDG 12}).

\section{Drone Technology and Operational Innovation} \label{sec:Innovation}

Drone technology has emerged as a catalyst for operational innovation, transforming the way businesses operate across various industries. Organizations strive for greater efficiency, productivity, and customer satisfaction, and drone technology has emerged as a disruptive force, revolutionizing operational processes. By leveraging their aerial capabilities, drones enable organizations to collect real-time data, automate tasks, and improve decision-making. As an edge example, Domino's collaborated with Gravity Industries to deliver pizza with a courier flying in a jetpack dress at the Glastonbury music festival in Somerset, United Kingdom \citep{heier2023}. In this section, we highlight the transformative power of drone technology and its role in shaping innovation through operations.

First, Drone-as-a-Service (DaaS) and subscription-based models are emerging, allowing organizations to access drone capabilities without the need for significant upfront investments in hardware, software, maintenance, or pilot training. For sectors that need drones for a temporary period (i.e., mine mapping), the renting option is beneficial \citep{kumas2022daas}. In contrast, for large-scale e-commerce companies (i.e., commercial last-mile delivery) or delivery services (i.e., commercial pickup and delivery delivery), investing in a drone fleet requires operational and financial planning. There needs to be an emphasis on how organizations will manage the change if they want to invest or rent a drone fleet. Because of the COVID-19 pandemic and the increase in remote working/studying, people started to move out of big cities and therefore the entire population is dispersed to a larger area \citep{levin2021insights}. With increased competition in e-commerce, companies should satisfy customers with minimized lead times. In order to sustain their same-day delivery and on-demand delivery, drones must be employed in their daily operations. Similarly, food delivery and the secondhand market are the most dominant commercial pickup and delivery operations. Many e-commerce companies are investing in second-hand trade platforms \citep{angueira2022}. With the increase in social media usage, platforms like Facebook marketplace and auction websites direct trade between people are becoming more popular. One major problem in the second-hand market is person-to-person shipping. Enabling a drone fleet for this direct trade could potentially create a competitive advantage. Drones can act as a tertiary party, like Uber Eats, in food delivery. A shared fleet of drones can serve multiple companies, for example, restaurants and grocery stores in a city. This would result in ride-sharing, assortment optimization, and pricing problems. If the drone fleet is in synchronization with other vehicles then the problems would be much more complex. People with their own individual drones can join the fleet under the gig economy model and create another income source for individuals.

Drone charging stations can be another stream of business, with the current expansion of electric car charging stations, the location of drone charging stations can be integrated into the decision-making process. The drone charging stations can also be integrated into the top parts of street lamps or parcel delivery lockers and this might generate partnerships between drone companies and stakeholders that are responsible for the manufacturing and maintenance of street lights and parcel delivery lockers. The availability of drones can change the location decision of parcel delivery lockers. For further energy efficiency, some drones use forms of public transportation (i.e., they land and dispatch from public buses and public trains) to handle freight requests. The scheduled line serves both passengers and drones, favoring passengers in case of capacity constraints \citep{mourad2021integrating}. This would create new planning directions and revenue streams for public transportation companies. They might even invest in their own drone fleet and get involved in e-commerce or delivery activities. However, with the increased use of drones, air traffic will be a new factor that needs to be considered in planning. The drones are making distance planning easier since they travel in the Euclidean distance, and generally, the distances are symmetric, enabling computational easiness. 

\section{Conclusion and Future Research Directions} \label{sec:Conclusion}

In conclusion, we review the academic literature on drone operations from 2017 to 2023 and classify them according to different industry applications (i.e., commercial last-mile delivery, commercial pickup and delivery, telecommunication, insurance, healthcare, humanitarian, environmental, urban planning, and homeland security). We analyze the articles in different industries in terms of commodity, network type (i.e., vehicle-drone synchronization or independent drone system), model, uncertainty, solution method, decision levels, and objective function. We then formally conceptualize the effect of drones on the United Nations Sustainable Development Goals and operational innovation. In short, organizations in different industries are employing drones to make their operations more sustainable and innovative for competitive advantage. The mechanical advancements in drones are making drones more autonomous. Therefore they are becoming less costly to manufacture and operate. The developing landscape of drone mechanics will continue to bolster and disrupt traditional operations, by creating alternatives and business models.

In terms of modeling, the current state of the art is highly focused on the computing aspect of the problem. Around 75\% of the articles use mixed integer linear programming, around 84\% of the articles consider deterministic settings and around 79\% of the articles use synthetic instances, instead of real-life data sets. The majority of the OM/OR/MS drone literature focuses on finding the optimal solution for large instances (up to 50 nodes) in a reasonable amount of time. In particular, commercial last-mile delivery operations and commercial pickup and delivery operations lack industry know-how and use more tools. Although healthcare and humanitarian drone operations articles were better at considering real-life case studies, there is a huge gap in equity measures when drones are involved in operations. Those classes of works used traditional minimized cost or minimized time objective functions and disregarded the equity dimension of the problems. Almost all of the articles consider single depot, which is an unrealistic assumption. The literature lacks analytical works that collaborate with actual stakeholders and works that use tools such as econometrics, field experiments, game theory, optimal control, and utility functions.

In addition, drone operations lack research in OM interfaces (i.e., OM-Finance, OM-Organizational Behaviour, OM-Marketing). For example, all of the works assume the drone fleet exists and do not consider the procurement process of drones. The investment phase is a good research area for the OM-Finance interface. How to price shared drones and how to handle assortment optimization are some other unanswered questions. On the other hand, for the OM-Organizational Behaviour interface, the effects of drones on the existing workforce or change management can be explored. 

Moreover, certain industries lack any analytical efforts concerning drone operations, such as but not limited to agriculture (i.e., crop management, pest and disease detection, irrigation, and fertilization optimization), construction (i.e., construction progress monitoring, safety inspections), green energy (i.e., inspection and maintenance of solar panels and wind turbines), mining (i.e., mining site and mineral detection) and telecommunication (i.e., measure signal quality). While tackling these industries, industry knowledge should be incorporated into drone operations for real-life applications. For instance, while routing the drone on land, the current humidity level of the land is a parameter that would affect the decision-making process. Finally, the effect of drones on supply chain financing, supply chain risk, and supply chain resilience can be further studied.

\ACKNOWLEDGMENT{}

\bibliographystyle{chicago}
\bibliography{references}

\ECSwitch

\section{Appendix} \label{sec:appendix}

\begin{table}[!hbt] \setlength{\tabcolsep}{8pt}
\begin{center}
\caption{Distribution of drone operations articles according to journals and industries.}
\begin{adjustbox}{width=0.8\textwidth}
\renewcommand{\arraystretch}{1.25}
\begin{tabular}{|p{3.75in}|p{0.15in}|p{0.15in}|p{0.15in}|p{0.15in}|p{0.15in}|p{0.15in}|p{0.15in}|p{0.15in}|p{0.15in}|p{0.15in}|p{0.15in}|}
\hline
\rowcolor{Gray}
& \multicolumn{10}{c}{\textbf{Industry}} & \\
\hline
\rowcolor{Gray}
\rotatebox{0}{\textbf{Journal Name}} & \rotatebox{270}{\textbf{Commercial Last-mile Delivery Operations}} & \rotatebox{270}{\textbf{Commercial Pickup and Delivery Operations     }} & \rotatebox{270}{\textbf{Telecommunications Operations}} & \rotatebox{270}{\textbf{Insurance Operations}} & \rotatebox{270}{\textbf{Healthcare Operations}} & \rotatebox{270}{\textbf{Humanitarian Operations}} & \rotatebox{270}{\textbf{Environmental Operations}} & \rotatebox{270}{\textbf{Urban Planning Operations}} & \rotatebox{270}{\textbf{Homeland Security Operations}} & \rotatebox{270}{\textbf{Methodological}} & \rotatebox{270}{\textbf{Journal Total}} \\
\hline
Transportation Research Part C & 9 & 5 &  & 1 & 1 & 2 & 1 &  &  & 9 & 28 \\
Computers and Operations Research & 6 & 2 & 1 &  & 2 &  & 1 & 2 &  & 2 & 16 \\
European Journal of Operational Research & 9 & 1 &  &  &  & 3 &  & 1 &  & 1 & 15 \\
Computers \& Industrial Engineering  & 6 & 2 &  &  & 4 & 1 &  &  & 1 &  & 14 \\
Transportation Science & 5 & 1 &  &  & 1 & 1 &  &  &  & 2 & 10 \\
Annals of Operations Research & 4 &  &  &  & 1 & 2 & 1 &  &  & 1 & 9 \\
Transportation Research Part E & 4 & 1 &  &  &  & 2 & 1 &  &  & 1 & 9 \\
Transportation Research Part B & 3 &  &  &  &  & 2 & 1 &  &  & 3 & 9 \\
Networks & 4 &  &  &  &  &  &  &  &  & 4 & 8 \\
OR Spectrum & 1 &  &  &  & 1 & 3 &  &  &  &  & 5 \\
International Journal of Production Economics & 2 &  &  &  & 1 & 1 &  &  &  &  & 4 \\
Omega &  & 1 &  &  & 1 &  &  &  &  & 1 & 3 \\
INFORMS Journal on Computing &  &  &  &  &  &  &  &  &  & 2 & 2 \\
Journal of Intelligent \& Robotic Systems &  &  &  &  & 2 &  &  &  &  &  & 2 \\
Manufacturing \& Service Operations Management & 1 &  &  &  & 1 &  &  &  &  &  & 2 \\
Optimization Letters & 2 &  &  &  &  &  &  &  &  &  & 2 \\
Production and Operations Management & 1 &  &  &  &  & 1 &  &  &  &  & 2 \\
Algorithms &  &  &  &  &  &  &  &  &  & 1 & 1 \\
Applied Energy & 1 &  &  &  &  &  &  &  &  &  & 1 \\
Applied Mathematical Modelling &  &  &  &  &  &  &  &  &  & 1 & 1 \\
Expert Systems with Applications & 1 &  &  &  &  &  &  &  &  &  & 1 \\
IEEE Transactions on Transportation Electrification & 1 &  &  &  &  &  &  &  &  &  & 1 \\
International Journal of Disaster Risk Reduction &  &  &  &  &  & 1 &  &  &  &  & 1 \\
International Journal of Logistics Research and Applications &  & 1 &  &  &  &  &  &  &  &  & 1 \\
International Journal of Production Research & 1 &  &  &  &  &  &  &  &  &  & 1 \\
International Transactions in Operational Research &  &  &  &  &  &  &  &  &  & 1 & 1 \\
Journal of Global Optimization &  &  & 1 &  &  &  &  &  &  &  & 1 \\
Journal of the Operational Research Society & 1 &  &  &  &  &  &  &  &  &  & 1 \\
Physical Communication &  &  &  &  &  &  &  &  &  & 1 & 1 \\
Sustainability & 1 &  &  &  &  &  &  &  &  &  & 1 \\
Transportation Research Part D &  &  &  &  &  &  & 1 &  &  &  & 1 \\
Vaccine &  &  &  &  & 1 &  &  &  &  &  & 1 \\
Preprint & 1 &  &  &  & 2 &  &  &  &  &  & 3 \\
\hline
\textbf{Industry Total} & 64 & 14 & 2 & 1 & 18 & 19 & 6 & 3 & 1 & 30 & 158 \\
\hline
\end{tabular}
\end{adjustbox}
\label{table:journals}
\end{center}
\end{table}

\begin{figure}
\begin{center}
\caption{Keyword analysis of commercial last-mile delivery operations articles ($n=64$).} \medskip
\includegraphics[width=0.8\textwidth]{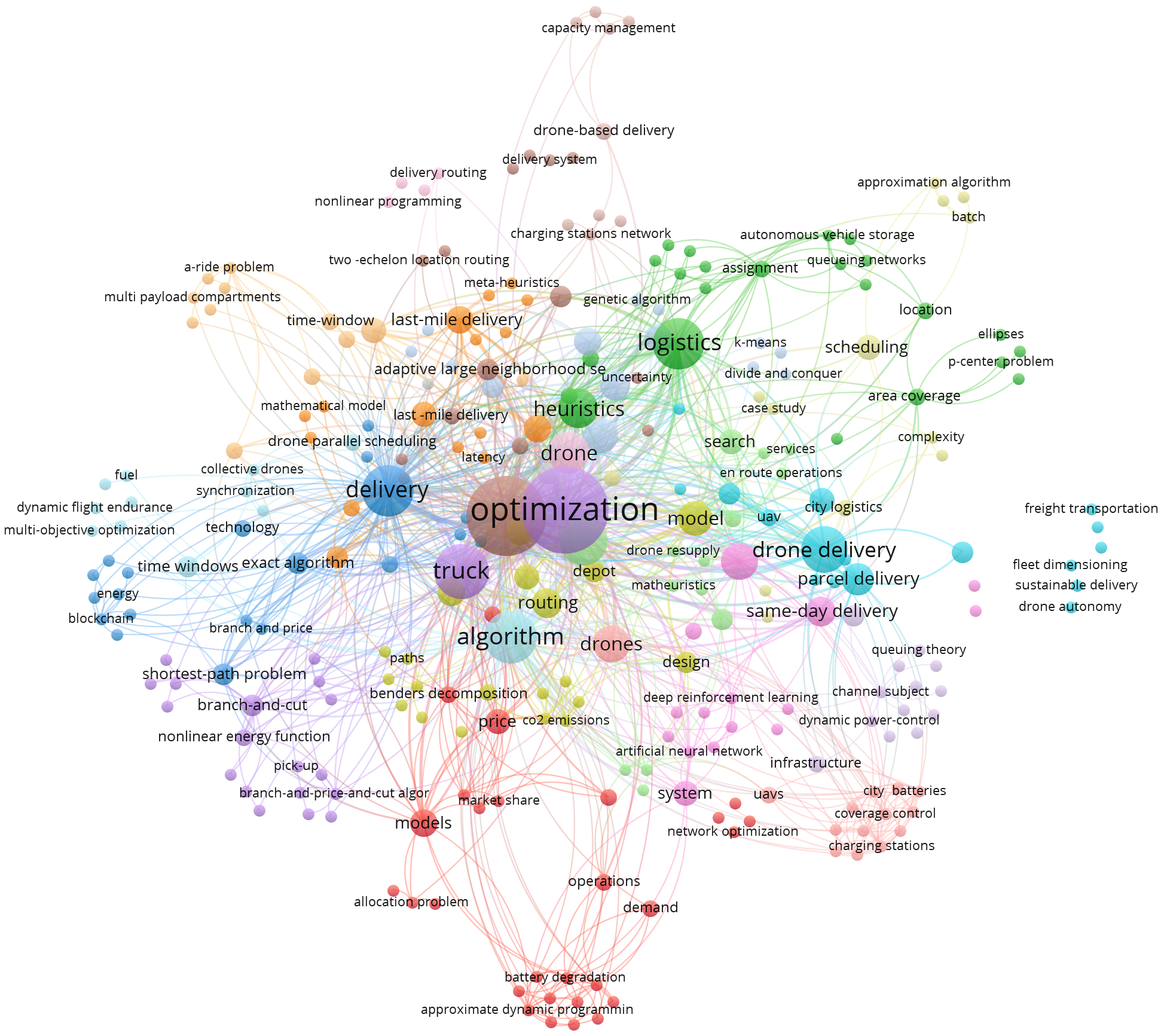}
\label{Img:last_mile_keywords}
\end{center}
\end{figure}

\begin{figure}
\begin{center}
\caption{Keyword analysis of commercial pickup and delivery operations articles ($n=14$).} \medskip
\includegraphics[width=0.8\textwidth]{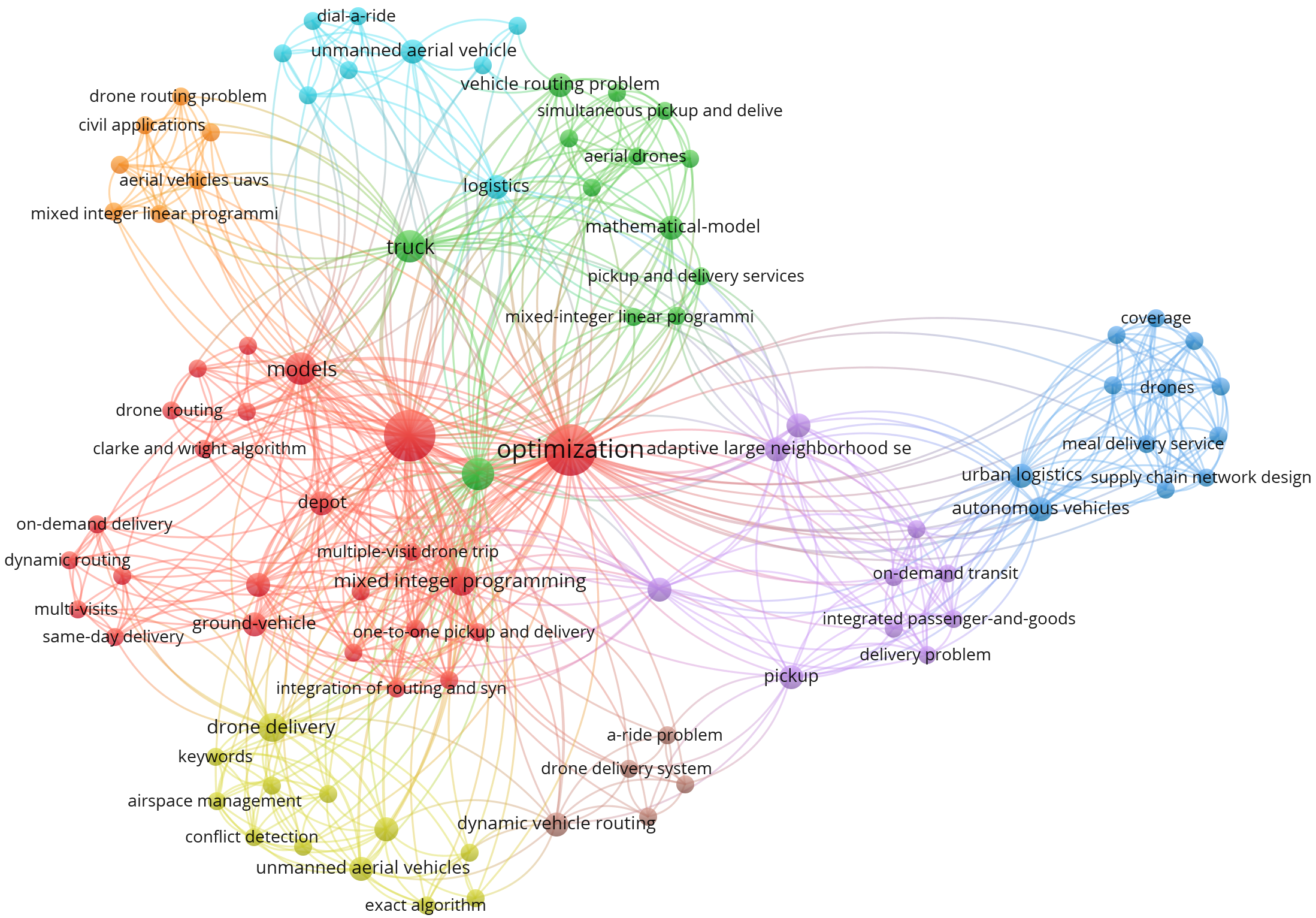}
\label{Img:pickup_keywords}
\end{center}
\end{figure}

\begin{figure}
\begin{center}
\caption{Keyword analysis of healthcare operations articles ($n=18$).} \medskip
\includegraphics[width=0.9\textwidth]{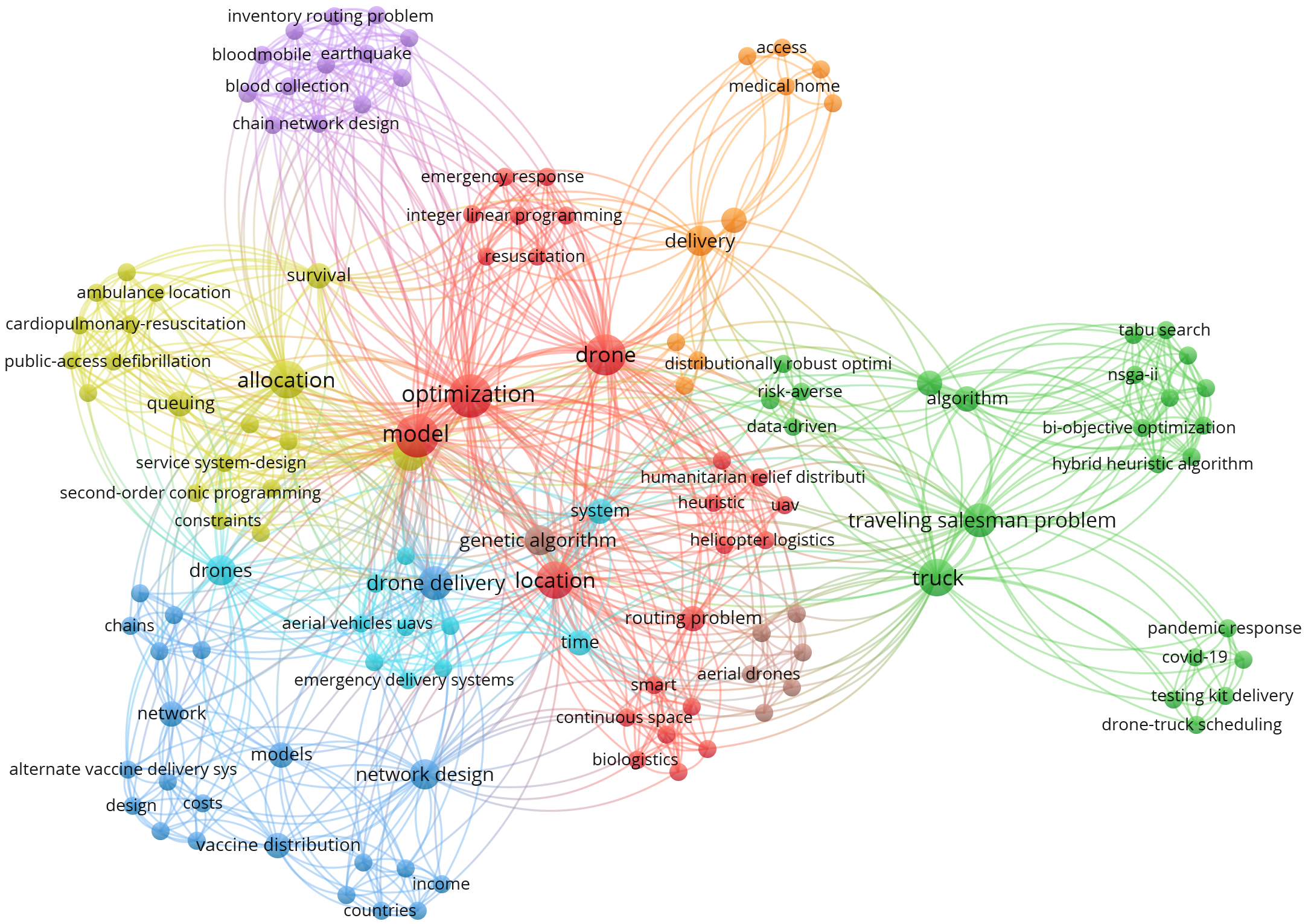}
\label{Img:healthcare_keywords}
\end{center}
\end{figure}

\begin{figure}
\begin{center}
\caption{Keyword analysis of humanitarian operations articles ($n=19$).} \medskip
\includegraphics[width=0.9\textwidth]{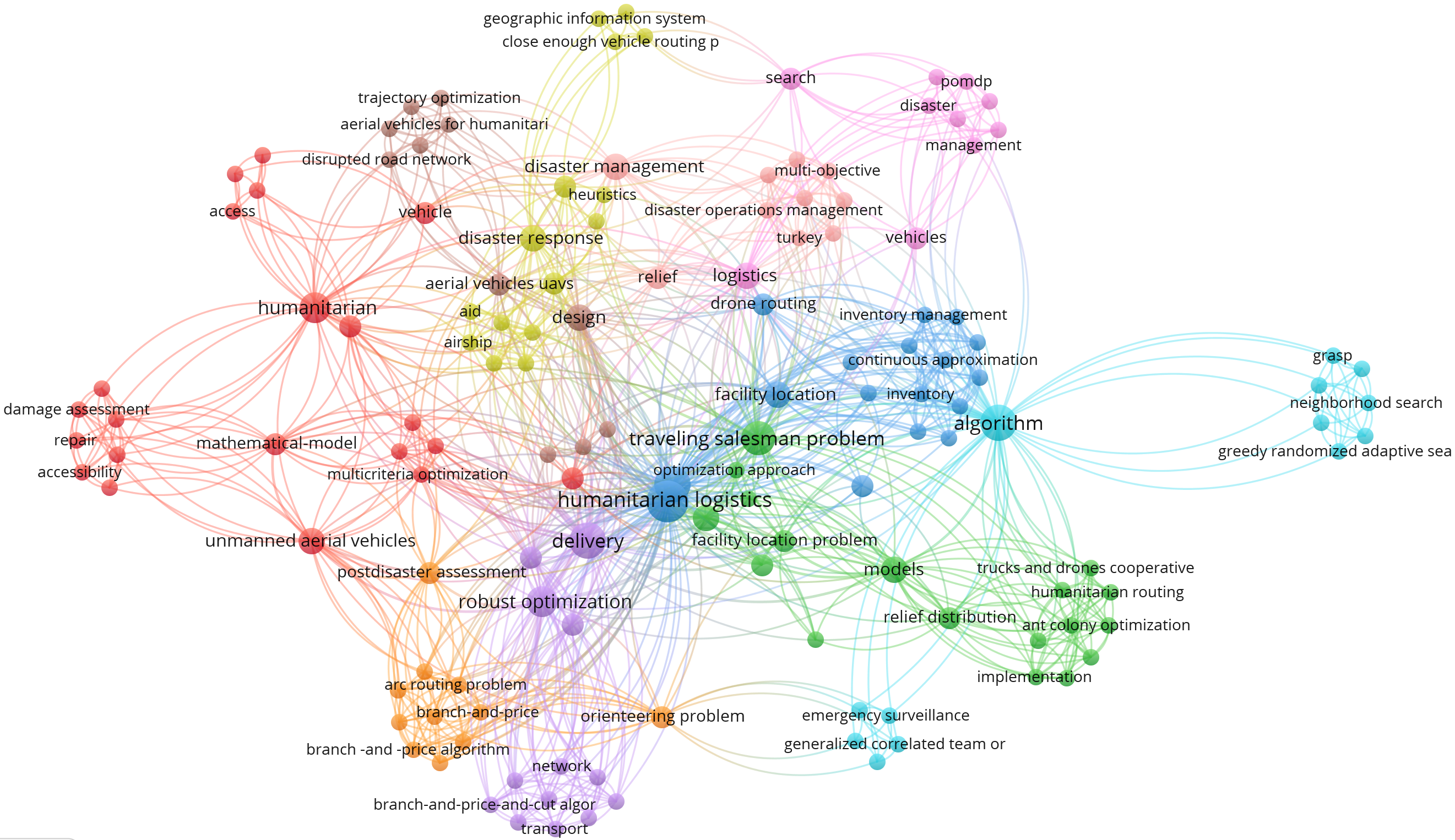}
\label{Img:humanitarian_keywords}
\end{center}
\end{figure}

\begin{figure}
\begin{center}
\caption{Keyword analysis of environmental operations articles ($n=6$).} \medskip
\includegraphics[width=0.9\textwidth]{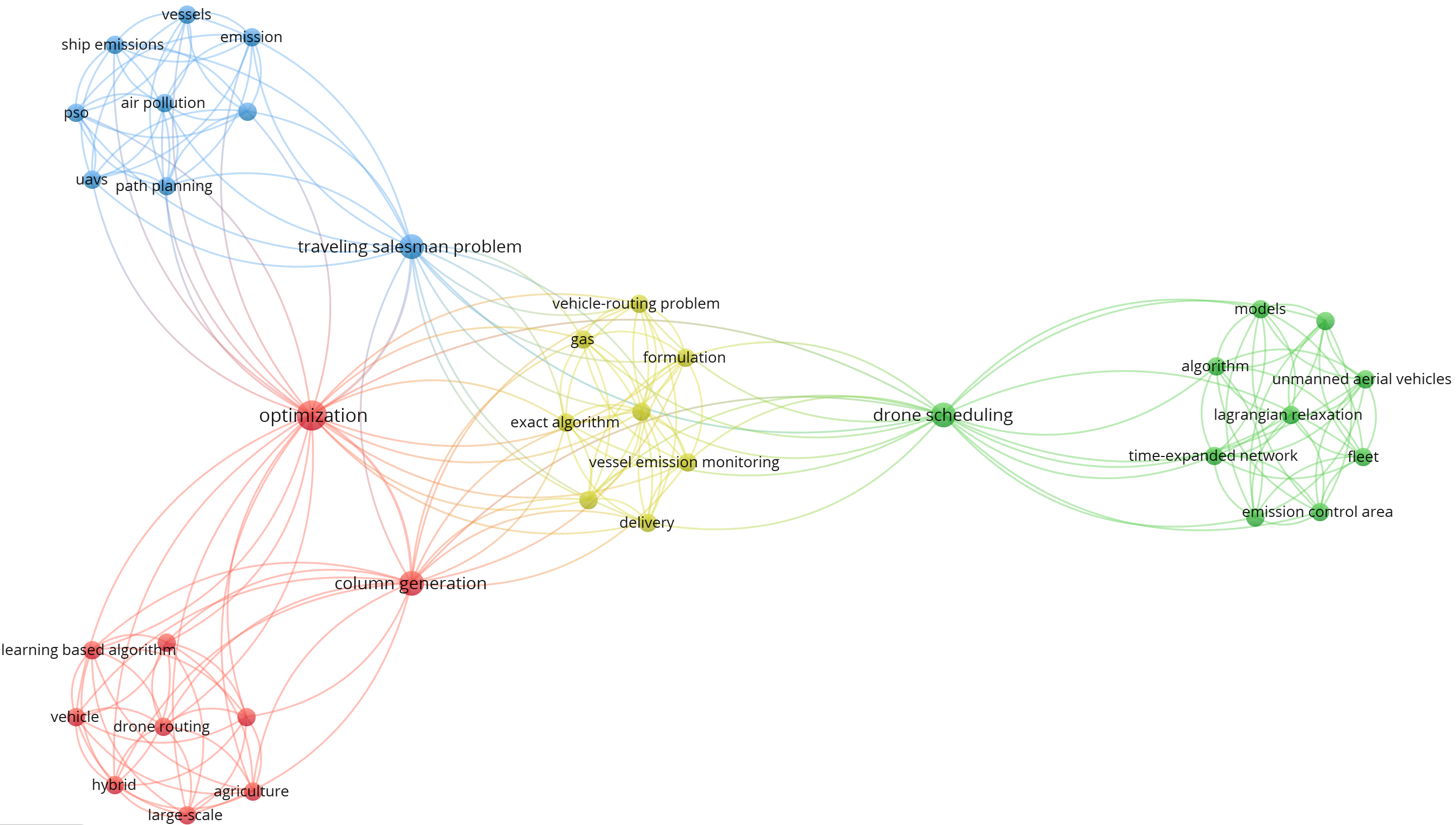}
\label{Img:environmental_keywords}
\end{center}
\end{figure}

\begin{figure}
\begin{center}
\caption{Keyword analysis of telecommunications, insurance, urban planning, and homeland security operations articles ($n=7$).} \medskip
\includegraphics[width=0.9\textwidth]{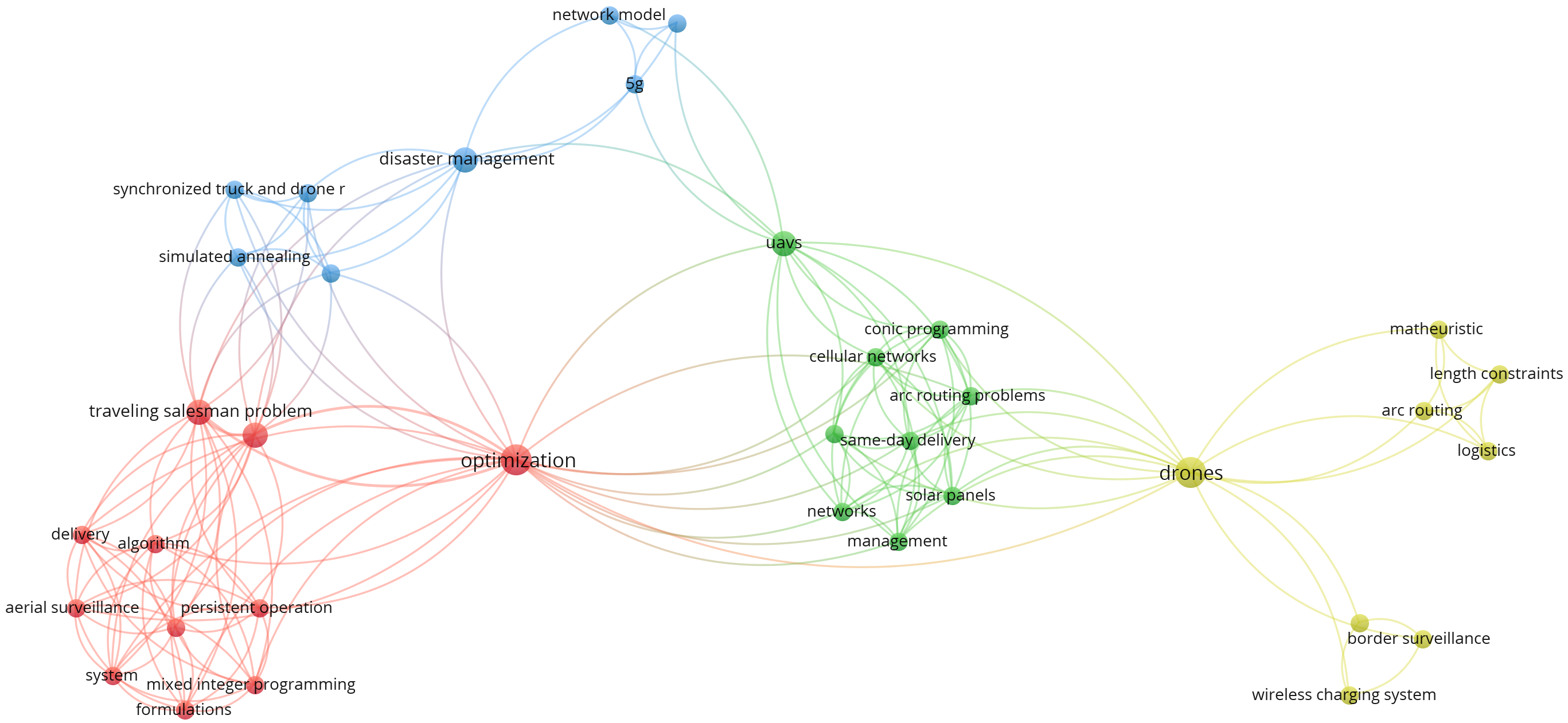}
\label{Img:other_keywords}
\end{center}
\end{figure}

\begin{figure}
\begin{center}
\caption{Keyword analysis of methodological articles ($n=30$).} \medskip
\includegraphics[width=0.8\textwidth]{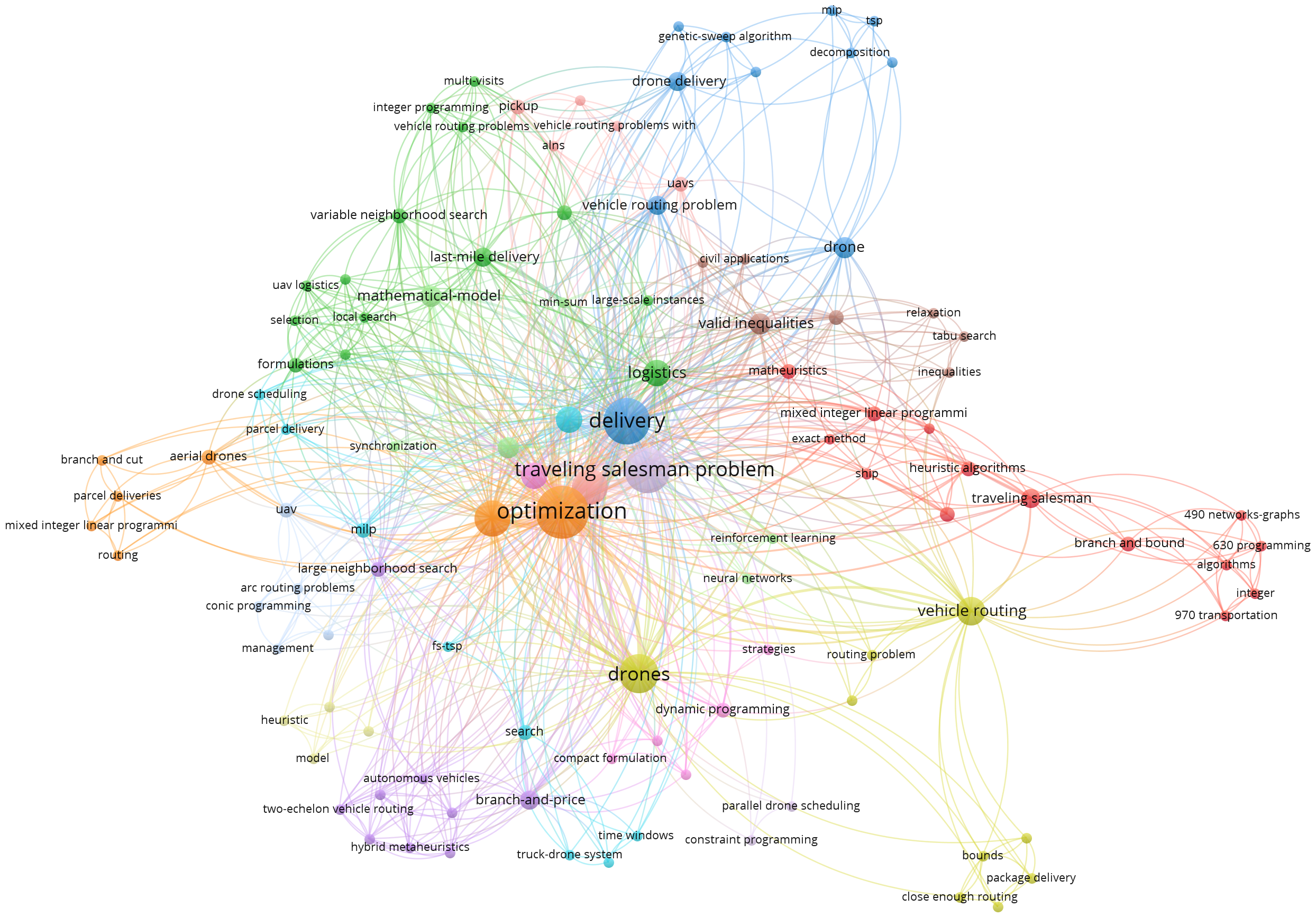}
\label{Img:methods_keywords}
\end{center}
\end{figure}

\begin{table}[htp] \setlength{\tabcolsep}{8pt}
\begin{center}
\caption{List of abbreviations used in Appendix.}
\begin{adjustbox}{width=\textwidth}
\begin{normalsize}
\renewcommand{\arraystretch}{1}
\begin{tabular}{llll}
\hline
\rowcolor{Gray}
\multicolumn{4}{l}{\textbf{Table \ref{table:detailed-features}}}\\
\hline
\multicolumn{4}{l}{\textbf{Commodity}}\\
\hline
\textit{Unit:} & AED: Automated External Defibrillator & B: Blood & DS: Disinfectant Spray \\
 & GHG: General Humanitarian Goods & GMG: General Medical Goods & Me: Medication \\
 & Meal: Meal & Pa: Package & TK: Testing Kit \\
 & V: Vaccine & Wi-Fi: Bandwidth &  \\
 \hdashline
\textit{Characteristic:} & (M): Multiple & (P): Perishable &  \\
\hline
\multicolumn{4}{l}{\textbf{Pickup}}\\
\hline
\end{tabular}
\end{normalsize}
\end{adjustbox}
\label{table:abbreviations}
\end{center}
\end{table}

\begin{table}[!hbt] \setlength{\tabcolsep}{8pt}
\begin{center}
\caption{Detailed features of drone operations articles.}
\begin{adjustbox}{width=\textwidth}
\renewcommand{\arraystretch}{1.25}
\begin{tabular}{|p{2in}|p{0.75in}|p{0.5in}|p{0.5in}|p{0.5in}|p{0.5in}|p{0.5in}|p{0.5in}|p{0.5in}|p{0.5in}|p{3in}|}
\hline
\rowcolor{Gray}
\rotatebox{0}{\textbf{Article}} & \rotatebox{270}{\textbf{Commodity}} & \rotatebox{270}{\textbf{Pickup}} & \rotatebox{270}{\textbf{Loop}} & \rotatebox{270}{\textbf{Synchronization Vehicle     }} & \rotatebox{270}{\textbf{Synchronization Type}} & \rotatebox{270}{\textbf{Recharging Node}} & \rotatebox{270}{\textbf{Number of Vehicles}} & \rotatebox{270}{\textbf{Number of Drones}} & \rotatebox{270}{\textbf{Nodes per Drone Trip}} & \rotatebox{0}{\textbf{Case Study}} \\
\hline
\multicolumn{11}{|l|}{\textbf{Commercial Last-mile Delivery Operations}} \\
\hline
\citet{asadi2022monotone} & Pa & $\times$ & $\times$ & $\times$ & $\times$ & $\checkmark$ & $\times$ & $\infty$ & $\times$ & Synthetic \\
\citet{bakir2020optimizing} & Pa & $\times$ & $\checkmark$ & T & Co & $\times$ & $> 1$ & $> 1$ & $1$ & Synthetic \\
\citet{baloch2020strategic} & Pa & $\times$ & $\odot$ & $\times$ & $\times$ & $\times$ & $\times$ & $> 1$ & $\times$ & New York City \\
\citet{betti2023dispatching} & Pa & $\times$ & $\times$ & $\times$ & $\times$ & $\times$ & $\times$ & $> 1$ & $\times$ & Synthetic \\
\citet{boysen2018drone} & Pa & $\times$ & $\times$ & T & Co & $\times$ & $1$ & $> 1$ & $1$ & Synthetic \\
\citet{bruni2022logic} & Pa & $\times$ & $\odot$ & T & Co & $\times$ & $1$ & $> 1$ & $1$ & Synthetic \\
\hline
\end{tabular}
\end{adjustbox}
\label{table:detailed-features}
\end{center}
\end{table}

\begin{table}[!hbt] \setlength{\tabcolsep}{8pt}
\begin{center}
\begin{adjustbox}{width=\textwidth}
\renewcommand{\arraystretch}{1.25}
\begin{tabular}{|p{2in}|p{0.75in}|p{0.5in}|p{0.5in}|p{0.5in}|p{0.5in}|p{0.5in}|p{0.5in}|p{0.5in}|p{0.5in}|p{3in}|}
\hline
\rowcolor{Gray}
\rotatebox{0}{\textbf{Article}} & \rotatebox{270}{\textbf{Commodity}} & \rotatebox{270}{\textbf{Pickup}} & \rotatebox{270}{\textbf{Loop}} & \rotatebox{270}{\textbf{Synchronization Vehicle     }} & \rotatebox{270}{\textbf{Synchronization Type}} & \rotatebox{270}{\textbf{Recharging Node}} & \rotatebox{270}{\textbf{Number of Vehicles}} & \rotatebox{270}{\textbf{Number of Drones}} & \rotatebox{270}{\textbf{Nodes per Drone Trip}} & \rotatebox{0}{\textbf{Case Study}} \\
\hline
\multicolumn{11}{|l|}{\textbf{Commercial Pickup and Delivery Operations}} \\
\hline
\citet{cheng2023adaptive} & Pa & $\checkmark$ & $\times$ & Bus & Se, Co & $\times$ & $> 1$ & $> 1$ & $1$ & Synthetic \\
\citet{gomez2022pickup} & Pa & $\checkmark$ & $\times$ & $\times$ & $\times$ & $\times$ & $\times$ & $> 1$ & $\times$ & Synthetic \\
\citet{gu2023dynamic} & Pa & $\checkmark$ & $\times$ & T & Co & $\times$ & $1$ & $1$ & $> 1$ & Synthetic \\
\citet{ham2018integrated} & Pa & $\checkmark$ & $\checkmark$ & T & Se, Co & $\times$ & $> 1$ & $\geq 1$ & $> 1$ & Synthetic \\
\citet{jiang2023multi} & Pa & $\checkmark$ & $\times$ & T & Co & $\times$ & $> 1$ & $> 1$ & $> 1$ & Synthetic \\
\citet{karak2019hybrid} & Pa (M) & $\checkmark$ & $\times$ & T & Co & $\times$ & $1$ & $> 1$ & $> 1$ & Synthetic \\
\hline
\end{tabular}
\end{adjustbox}
\label{table:articles}
\end{center}
\end{table}

\begin{table}[!hbt] \setlength{\tabcolsep}{8pt}
\begin{center}
\begin{adjustbox}{width=\textwidth}
\renewcommand{\arraystretch}{1.25}
\begin{tabular}{|p{2in}|p{0.75in}|p{0.5in}|p{0.5in}|p{0.5in}|p{0.5in}|p{0.5in}|p{0.5in}|p{0.5in}|p{0.5in}|p{3in}|}
\hline
\rowcolor{Gray}
\rotatebox{0}{\textbf{Article}} & \rotatebox{270}{\textbf{Commodity}} & \rotatebox{270}{\textbf{Pickup}} & \rotatebox{270}{\textbf{Loop}} & \rotatebox{270}{\textbf{Synchronization Vehicle     }} & \rotatebox{270}{\textbf{Synchronization Type}} & \rotatebox{270}{\textbf{Recharging Node}} & \rotatebox{270}{\textbf{Number of Vehicles}} & \rotatebox{270}{\textbf{Number of Drones}} & \rotatebox{270}{\textbf{Nodes per Drone Trip}} & \rotatebox{0}{\textbf{Case Study}} \\
\hline
\multicolumn{11}{|l|}{\textbf{Urban Planning Operations}} \\
\hline
\citet{amorosi2023multiple} & $\times$ & $\times$ & $\times$ & M & Co & $\times$ & $1$ & $> 1$ & $> 1$ & Cordoba Courtyard Festival \\
\citet{campbell2021solving} & $\times$ & $\times$ & $\times$ & $\times$ & $\times$ & $\times$ & $\times$ & $> 1$ & $\times$ & Synthetic \\
\citet{manshadian2023synchronized} & DS & $\times$ & $\checkmark$ & T & Co & $\times$ & $> 1$ & $> 1$ & $> 1$ & Tehran, Iran (Environmental Contamination) \\
\hline
\multicolumn{11}{|l|}{\textbf{Homeland Security Operations}} \\
\hline
\citet{ahmadian2022smart} & $\times$ & $\times$ & $\times$ & $\times$ & $\times$ & $\checkmark$ &  $\times$ & $> 1$ & $\times$ & US-Mexico Border \\
\hline
\multicolumn{11}{|l|}{\textbf{Methodological}} \\
\hline
\end{tabular}
\end{adjustbox}
\label{table:articles2}
\end{center}
\end{table}

\begin{table}[!hbt] \setlength{\tabcolsep}{8pt}
\begin{center}
\caption{Modeling, uncertainty and solution approaches considered in each drone operations article.}
\begin{adjustbox}{width=\textwidth}
\renewcommand{\arraystretch}{1.25}
\begin{tabular}{|p{2in}|p{1in}|p{0.6in}|p{6in}|}
\hline
\rowcolor{Gray}
\rotatebox{0}{\textbf{Article}} & \rotatebox{270}{\textbf{Model}} & \rotatebox{270}{\textbf{Uncertainty     }} & \rotatebox{0}{\textbf{Solution Method}} \\
\hline
\multicolumn{4}{|l|}{\textbf{Commercial Last-mile Delivery Operations}} \\
\hline
\citet{asadi2022monotone} & FHMDP  & D & [H] Monotone approximate dynamic programming with regression-based initialization. \\
\citet{bakir2020optimizing} & MILP & $\times$ & [E] Dynamic discretization discovery approach. \\
\citet{baloch2020strategic} & MINLP, U & $\times$ & [E] Logic-based Benders decomposition. \\
\citet{betti2023dispatching} & RA & $\times$ & [H] Approximation algorithm. \\
\citet{boysen2018drone} & MILP & $\times$ & [E] Decomposition algorithm. [H] Simulated annealing. \\
\citet{bruni2022logic} & MILP & $\times$ & [E] Logic-based Benders decomposition. \\
\hline
\end{tabular}
\end{adjustbox}
\label{table:models-uncertainties-solutions}
\end{center}
\end{table}

\begin{table}[!hbt] \setlength{\tabcolsep}{8pt}
\begin{center}
\begin{adjustbox}{width=\textwidth}
\renewcommand{\arraystretch}{1.25}
\begin{tabular}{|p{2in}|p{1in}|p{0.6in}|p{6in}|}
\hline
\rowcolor{Gray}
\rotatebox{0}{\textbf{Article}} & \rotatebox{270}{\textbf{Model}} & \rotatebox{270}{\textbf{Uncertainty     }} & \rotatebox{0}{\textbf{Solution Method}} \\
\hline
\multicolumn{4}{|l|}{\textbf{Commercial Last-mile Delivery Operations (Cont'd)}} \\
\hline
\citet{el2021parcel} & MILP & $\times$ & [H] Valid inequalities. \\
\citet{huang2020method} & MILP  & $\times$ & [E] Full coverage deployment, local movement, and adjustment algorithms. \\
\citet{jeong2019truck} & MILP  & $\times$ & [H] Genetic algorithm, particle swarm optimization, simulated annealing, and nearest neighborhood.\\
\citet{jung2022drone} & RO & WD & [E] Reformulation, commercial solver. \\
\citet{kang2021exact} & MILP & $\times$ & [E] Logic-based Benders decomposition. \\
\citet{kitjacharoenchai2019multiple} & MILP & $\times$ & [H] Adaptive insertion algorithm. \\
\hline
\end{tabular}
\end{adjustbox}
\label{table:articles}
\end{center}
\end{table}

\begin{table}[!hbt] \setlength{\tabcolsep}{8pt}
\begin{center}
\begin{adjustbox}{width=\textwidth}
\renewcommand{\arraystretch}{1.25}
\begin{tabular}{|p{2in}|p{1in}|p{0.6in}|p{6in}|}
\hline
\rowcolor{Gray}
\rotatebox{0}{\textbf{Article}} & \rotatebox{270}{\textbf{Model}} & \rotatebox{270}{\textbf{Uncertainty     }} & \rotatebox{0}{\textbf{Solution Method}} \\
\hline
\multicolumn{4}{|l|}{\textbf{Telecommunications Operations}} \\
\hline
\citet{cicek2019location} & MINLP  & $\times$ & [H] Modified alternate location-allocation algorithm. \\
\citet{colajanni2023three} &  ThSSP  & D & [E] Variational formulation. \\ 
\hline
\multicolumn{4}{|l|}{\textbf{Insurance Operations}} \\
\hline
\citet{zeng2022nested} & MIQCP  & $\times$ & [H] MILP reformulation, neighborhood search algorithms. \\
\hline
\multicolumn{4}{|l|}{\textbf{Healthcare Operations}} \\
\hline
\citet{boutilier2022drone} &  MILP, Q  & RT  & [E] Reformulation, commercial solver. \\
\hline
\end{tabular}
\end{adjustbox}
\label{table:articles}
\end{center}
\end{table}

\begin{table}[!hbt] \setlength{\tabcolsep}{8pt}
\begin{center}
\begin{adjustbox}{width=\textwidth}
\renewcommand{\arraystretch}{1.25}
\begin{tabular}{|p{2in}|p{1in}|p{0.6in}|p{6in}|}
\hline
\rowcolor{Gray}
\rotatebox{0}{\textbf{Article}} & \rotatebox{270}{\textbf{Model}} & \rotatebox{270}{\textbf{Uncertainty     }} & \rotatebox{0}{\textbf{Solution Method}} \\
\hline
\multicolumn{4}{|l|}{\textbf{Methodological}} \\
\hline
\citet{agatz2018optimization} & MILP  & $\times$ & [H] Route-first, cluster-second algorithm based on local search and dynamic programming. \\
\citet{amorosi2021coordinating} & MISOCP  & $\times$ & [H] Decomposition based algorithm. \\
\citet{boccia2021column} & MILP & $\times$ & [E] Branch-and-cut, column generation, valid inequalities. \\
\citet{boccia2023new} & MILP & $\times$ & [E] Arc-based integer linear programming reformulation, branch-and-cut algorithm. \\
\citet{bogyrbayeva2023deep} & MILP & $\times$ & [H] Deep reinforcement learning and neural networks. \\
\citet{bouman2018dynamic} & MILP & $\times$ & [E] Dynamic programming. \\
\hline
\end{tabular}
\end{adjustbox}
\label{table:articles}
\end{center}
\end{table}

\begin{table}[!hbt] \setlength{\tabcolsep}{8pt}
\begin{center}
\caption{Decision levels and objective functions considered in each drone operations article.}
\begin{adjustbox}{width=\textwidth}
\renewcommand{\arraystretch}{1.25}
\begin{tabular}{|p{2in}|p{0.15in}|p{0.15in}|p{0.15in}|p{0.15in}|p{0.15in}|p{3in}|p{2in}|p{0.15in}|p{0.15in}|p{0.15in}|p{0.15in}|p{0.15in}|p{3in}|}
\hline
\rowcolor{Gray}
\rotatebox{0}{\textbf{Article}} & \rotatebox{270}{\textbf{Location}} & \rotatebox{270}{\textbf{Inventory}} & \rotatebox{270}{\textbf{Allocation}} & \rotatebox{270}{\textbf{Distribution}} & \rotatebox{270}{\textbf{Routing}} & \rotatebox{270}{\textbf{Objective Function     }} & \rotatebox{0}{\textbf{Article}} & \rotatebox{270}{\textbf{Location}} & \rotatebox{270}{\textbf{Inventory}} & \rotatebox{270}{\textbf{Allocation}} & \rotatebox{270}{\textbf{Distribution}} & \rotatebox{270}{\textbf{Routing}} & \rotatebox{270}{\textbf{Objective Function     }} \\
\hline
\multicolumn{7}{|l|}{\textbf{Commercial Last-mile Delivery Operations}} & \multicolumn{7}{|l|}{\textbf{Insurance Operations}} \\ 
\hline
\citet{asadi2022monotone} & $\times$  & $\checkmark$  & $\checkmark$  & $\checkmark$  & $\times$  & [Max] p & \citet{zeng2022nested} & $\checkmark$  & $\times$  & $\checkmark$  & $\times$  & $\checkmark$  & [Min] ms \\ 
\cline{8-14}
\citet{bakir2020optimizing} & $\times$  & $\times$  & $\times$  & $\checkmark$  & $\checkmark$  & [Min] ms & \multicolumn{7}{|l|}{\textbf{Healthcare Operations}} \\ 
\cline{8-14}
\citet{baloch2020strategic} & $\checkmark$  & $\times$  & $\checkmark$  & $\times$  & $\times$  & [Max] p $^{M1,M2,M3,M4}$ & \citet{boutilier2022drone} &  $\checkmark$  & $\checkmark$  & $\checkmark$  & $\times$  & $\times$  & [Min] nd \\ 
\citet{betti2023dispatching} & $\checkmark$  & $\times$  & $\times$  & $\times$  & $\times$  & [Min] ds & \citet{dhote2020designing} &  $\checkmark$  & $\times$  & $\checkmark$  & $\times$  & $\times$  & [Min] tv $^{M1,M2,M3,M4}$ \\ 
\citet{boysen2018drone} & $\times$  & $\times$  & $\times$  & $\times$  & $\checkmark$  & [Min] ms $^{M1,M2}$ & \citet{enayati2023multimodal} &  $\checkmark$  & $\times$  & $\checkmark$  & $\checkmark$  & $\times$  & [Min] tc \\
\hline
\end{tabular}
\end{adjustbox}
\label{table:DL-obj}
\end{center}
\end{table}

\FloatBarrier

\begin{table}[h] \setlength{\tabcolsep}{2pt}
\caption{United Nations Sustainable Development Goals}
\begin{center}
\begin{normalsize}
\renewcommand{\arraystretch}{1.25}
\begin{tabular}{cll}
\hline
\textbf{UNSDG}&\textbf{Description}&\textbf{Class}\\
\hline
SDG 1 & No Poverty & Social \\
SDG 2 & Zero Hunger & Social \\
SDG 3 & Good Health and Well-being & Social \\
SDG 4 & Quality Education & Social \\
SDG 5 & Gender Equality & Social \\
SDG 6 & Clean Water and Sanitation & Environmental \\
SDG 7 & Affordable and Clean Energy & Environmental \\
SDG 8 & Decent Work and Economic Growth & Economic \\
SDG 9 & Industry, Innovation and Infrastructure & Economic \\
SDG 10 & Reduced Inequalities & Social \\
SDG 11 & Sustainable Cities and Communities & Social \\
SDG 12 & Responsible Consumption and Production & Economic \\
SDG 13 & Climate Action & Environmental \\
SDG 14 & Life Below Water & Environmental \\
SDG 15 & Life on Land & Environmental \\
SDG 16 & Peace, Justice and Strong Institutions & Social \\
SDG 17 & Partnerships & Social \\
\hline
\end{tabular}
\end{normalsize}
\end{center}
\label{table:UNSDG}
\end{table}

\end{document}